\documentclass[12pt]{amsart}%
\usepackage{amsmath}
\usepackage{amsfonts}
\usepackage{amssymb}
\usepackage{graphicx}
\usepackage{color}%
\providecommand{\U}[1]{\protect\rule{.1in}{.1in}}
\def\U{\Upsilon}

\newtheorem{theo}{Theorem}[section]
\newtheorem{prop}[theo]{Proposition}

\newtheorem{lem}[theo]{Lemma}

\begin{document}
\title[Degenerate reaction-diffusion system]{Exponential decay toward equilibrium \\via log convexity for a degenerate reaction-diffusion system}
\author{Laurent Desvillettes}
\address{Universit\'{e} de Paris and Sorbonne Universit\'e, CNRS, Institut de
Math\'{e}matiques de Jussieu-Paris Rive Gauche, F-75013, Paris, France}
\email{desvillettes@math.univ-paris-diderot.fr}
\author{Kim Dang Phung}
\address{Institut Denis Poisson, Universit\'{e} d'Orl\'{e}ans, Universit\'{e} de Tours
\& CNRS UMR 7013, B\^{a}timent de Math\'{e}matiques, Rue de Chartres, BP.
6759, 45067 Orl\'{e}ans, France}
\email{kim\_dang\_phung@yahoo.fr}
\date{}
\maketitle

\begin{abstract}
We consider a system of two reaction-diffusion equations coming out of
reversible chemistry. When the reaction happens on the totality of the domain,
it is known that exponential convergence to equilibrium holds. We show in this
paper that this exponential convergence also holds when the reaction holds
only on a given open set of a ball, thanks to an observation estimate deduced
by logarithmic convexity.

\end{abstract}

\bigskip

\bigskip

\section{Introduction and main result}

\bigskip

We are interested in a reaction-diffusion system coming out of the reversible
chemical reaction
\[
{\mathcal{A}}+{\mathcal{A}}\qquad\iff\qquad{\mathcal{B}}+{\mathcal{B}}\text{
,}%
\]
where ${\mathcal{A}}$ and ${\mathcal{B}}$ are chemical species of respective
concentrations $a:=a(x,t)\geq0$, $b:=b(x,t)\geq0$. \medskip

We suppose that the species ${\mathcal{A}}$ has a diffusion rate $d_{1}>0$ and
that the species ${\mathcal{B}}$ has a diffusion rate $d_{2}>0$. We also
assume that those species are confined in a chemical reactor represented by
the ball $\Omega:=B\left(  0,R\right)  :=\left\{  x\in\mathbb{R}^{\text{n}%
};\left\vert x\right\vert <R\right\}  $, where n$\in\{1,2,3\}$, and
$|\Omega|=1$ (that is $R=\frac{1}{2}$ if n$=1$, $R=\pi^{-1/2}$ if n$=2$, and
$R=(\frac{3}{4\pi})^{1/3}$ if n$=3$), so that homogeneous Neumann boundary
conditions are imposed. Finally and most importantly, the terms arising from
the reaction process are given by the mass action law, and the reaction rate
is proportional to the concentration of a catalyst $k:=k(x,t)\geq0$. This
paper is devoted to the treatment of the case when $k$ is not strictly
positive on $\Omega$. The system writes \medskip%

\begin{equation}
\left\{
\begin{array}
[c]{ll}%
{\partial}_{t}a-d_{1}\Delta a=k\left(  b^{2}-a^{2}\right)  \text{ ,} &
\quad\text{in}~\Omega\times\left(  0,+\infty\right)  \text{ ,}\\
{\partial}_{t}b-d_{2}\Delta b=-k\left(  b^{2}-a^{2}\right)  \text{ ,} &
\quad\text{in}~\Omega\times\left(  0,+\infty\right)  \text{ ,}\\
\partial_{n}a=\partial_{n}b=0\text{ ,} & \quad\text{on}~\partial\Omega
\times\left(  0,+\infty\right)  \text{ ,}\\
a\left(  \cdot,0\right)  =a_{0}\text{ , }b\left(  \cdot,0\right)  =b_{0}\text{
,} & \quad\text{in}~\Omega\text{ ,}%
\end{array}
\right.  \label{sys1}%
\end{equation}
where $\partial_{n}:=n(x)\cdot\nabla$, and $n(x)$ is the unit outward normal
vector at point $x\in\partial\Omega$. \medskip

We consider initial data $a_{0},b_{0}\in C^{2}\left(  \overline{\Omega
}\right)  $ (compatible with the Neumann boundary condition) which satisfy the
bound:
\begin{equation}
\forall x\in\overline{\Omega}\text{ ,}\quad0<B_{0}\leq a_{0}\left(  x\right)
\quad\text{and}\quad0<B_{0}\leq b_{0}\left(  x\right)  \text{ ,}
\label{hypin1}%
\end{equation}
for some constant $B_{0}>0$, and we suppose that
\begin{equation}
\int_{\Omega}\left[  a_{0}\left(  x\right)  +b_{0}\left(  x\right)  \right]
dx=2\text{ .} \label{hypin2}%
\end{equation}

At this point, we remark that at the formal level, the following \textit{{a
priori}} estimates hold:
\begin{equation}
\frac{d}{dt}\int_{\Omega}(a+b)=0\text{ ,} \label{ape1}%
\end{equation}

\begin{equation}
\label{ape2}\frac{1}{2}\frac{d}{dt}\left\Vert \left(  a, b\right)  \right\Vert
_{\left(  L^{2}\left(  \Omega\right)  \right)  ^{2}}^{2}+d_{1}\int_{\Omega
}\left\vert \nabla a\right\vert ^{2}+d_{2}\int_{\Omega}\left\vert \nabla
b\right\vert ^{2}+\int_{\Omega}k\left(  a+b\right)  \left\vert a -
b\right\vert ^{2}=0\text{ .}%
\end{equation}

Because of the terms $d_{1}\int_{\Omega}\left\vert \nabla a\right\vert ^{2}$
and $d_{2}\int_{\Omega}\left\vert \nabla b\right\vert ^{2}$, we expect that
$\lim_{t\rightarrow\infty}a(t,x)=a_{\infty}$, $\lim_{t\rightarrow\infty
}b(t,x)=b_{\infty}$, for some constants $a_{\infty}\geq0$ and $b_{\infty}%
\geq0$. Moreover, as soon as $k$ is strictly positive on some open set of
$\Omega$, $a_{\infty}=b_{\infty}$ because of the term $\int_{\Omega}k\left(
a+b\right)  \left\vert a-b\right\vert ^{2}$. Finally, estimate (\ref{ape1})
ensures that $\int_{\Omega}(a_{\infty}+b_{\infty})=\int_{\Omega}\left(
a_{0}+b_{0}\right)  =2$.

Remembering that $|\Omega|=1$, we finally expect that the equilibrium is
$\left(  a_{\infty},b_{\infty}\right)  =\left(  1,1\right)  $. \bigskip

It is proven thanks to an entropy method in \cite{DF} that when $k\left(
x,t\right)  \geq k_{0}>0$ for any $\left(  x,t\right)  \in\Omega\times\left(
0,+\infty\right)  $ (and for a domain $\Omega$ much more general than a ball),
then for any $t\geq0$,
\[
\left\Vert a\left(  \cdot,t\right)  -1\right\Vert _{L^{2}\left(
\Omega\right)  }^{2}+\left\Vert b\left(  \cdot,t\right)  -1\right\Vert
_{L^{2}\left(  \Omega\right)  }^{2}\leq\gamma\,e^{-\beta t}\left(  \left\Vert
a_{0}-1\right\Vert _{L^{2}\left(  \Omega\right)  }^{2}+\left\Vert
b_{0}-1\right\Vert _{L^{2}\left(  \Omega\right)  }^{2}\right)  \text{ ,}%
\]
where $\gamma>0$ and $\beta>0$ are constants which can be explicitly estimated
in terms of the data of the problem (that is $\Omega$, $k$, $d_{1}$, $d_{2}$).
\bigskip

Our main result extends this result in the case in which $k$ can be equal to
$0$ in part of a ball:

\begin{theo}
\label{t11} We define $\Omega:=B\left(  0,R\right)  :=\left\{  x\in
\mathbb{R}^{\normalfont{\text{n}}};\left\vert x\right\vert <R\right\}  $,
where $\normalfont{\text{n}}\in\{1,2,3\}$, the centered ball of $\mathbb{R}%
^{\normalfont{\text{n}}}$ of measure $1$. We also assume that $d_{1},d_{2}>0$
and $k\in C^{2}({\overline{\Omega}}\times\mathbb{R}_{+};\mathbb{R}_{+})$. We consider
initial data $a_{0},b_{0}\in C^{2}\left(  \overline{\Omega}\right)  $
(compatible with the Neumann boundary condition) which satisfy (\ref{hypin1})
and (\ref{hypin2}).

We finally assume that there exists $x_{0}\in\Omega$ and $r>0$ such that
$B\left(  x_{0},r\right)  \subset\Omega$ and $k\left(  x,t\right)  \geq
k_{0}>0$ for any $\left(  x,t\right)  \in B\left(  x_{0},r\right)
\times\left(  0,+\infty\right)  $. \medskip

Then there exists a unique smooth ($C^{2}({\overline{\Omega}}\times
\lbrack0,+\infty))$) solution to system (\ref{sys1}), such that
\begin{equation}
\inf_{t\geq0,x\in\overline{\Omega}}a(x,t)\geq B_{0}\text{ ,}\qquad\inf
_{t\geq0,x\in\overline{\Omega}}b(x,t)\geq B_{0}\text{ ,} \label{pmax}%
\end{equation}
and for any $t\geq0$,
\begin{equation}
\left\Vert a\left(  \cdot,t\right)  -1\right\Vert _{L^{2}\left(
\Omega\right)  }^{2}+\left\Vert b\left(  \cdot,t\right)  -1\right\Vert
_{L^{2}\left(  \Omega\right)  }^{2}\leq\gamma\,e^{-\beta t}\left(  \left\Vert
a_{0}-1\right\Vert _{L^{2}\left(  \Omega\right)  }^{2}+\left\Vert
b_{0}-1\right\Vert _{L^{2}\left(  \Omega\right)  }^{2}\right)  \text{ ,}
\label{mt}%
\end{equation}
where $\gamma>0$ and $\beta>0$ can be explicitly estimated in terms of
$|x_{0}|$, $r$, $\left\Vert \left(  a_{0},b_{0}\right)  \right\Vert
_{L^{3}\left(  \Omega\right)  }$, $B_{0}$, $k_{0}$, $||k||_{L^{\infty}\left(
\Omega\times\mathbb{R}_{+}\right)  }$, $d_{1}$, $d_{2}$.
\end{theo}

This theorem enables to explore a new issue, when one studies the convergence
towards an homogeneous chemical equilibrium for systems of reaction-diffusion
coming out of reversible chemistry. Indeed, the results of \cite{DF, DF3}
proven in the case of reaction-diffusion systems of 2, 3, or 4 equations with
strictly positive diffusion rates and reaction rates have been extended to a
large class of system in \cite{DFT}. The case when one diffusion rate is zero,
for simple systems, has also been studied in \cite{DF2} (see also \cite{BDS},
Rmk. 3.7), but the case of a degenerate reaction term had not been 
investigated, up to our knowledge.
\par
We wish to point out that the case of a system of two reaction-diffusion
equations in which one diffusion rate is zero and the reaction rate is also
zero on a non-negligeable set, is quite different, since the appearance of an
homogeneous equilibrium is not expected in such a situation.

\section{Proof of Theorem \ref{t11}}

We start here the \medskip

\textbf{{Proof of Theorem \ref{t11}}}: Note first that under the assumptions
of Theorem \ref{t11}, the existence, uniqueness and smoothness of a solution
to (\ref{sys1}) is a consequence of standard theorems for parabolic equations
(cf. for example \cite{D} or \cite{LSU}). The minimum principle estimate
(\ref{pmax}) can easily be seen at the formal level by considering, for a
given time $t$, the point $x \in\overline{\Omega}$ where $\min(a,b)$ reaches
its minimum. For a rigorous proof, we refer to \cite{D}, p.100-101 (where the
proof is detailed for a slightly different system). \medskip

Note then that for this solution, the following estimate holds:
\[
\frac{1}{3}\,\frac{d}{dt}\int_{\Omega}(a^{3}+b^{3})=-2d_{1}\int_{\Omega
}a\,|\nabla a|^{2}-2d_{2}\int_{\Omega}b\,|\nabla b|^{2}-\int_{\Omega}%
k\,(b^{2}-a^{2})^{2}\leq0\text{ .}%
\]
As a consequence%
\begin{equation}
\forall t\geq0\text{ ,}\qquad\int_{\Omega}\left[  a^{3}\left(  \cdot,t\right)
+b^{3}\left(  \cdot,t\right)  \right]  \leq\int_{\Omega}\left[  a_{0}%
^{3}\left(  x\right)  +b_{0}^{3}\left(  x\right)  \right]  dx\text{
.}\label{fp1}%
\end{equation}
\medskip

The proof of Theorem \ref{t11} is then an application of the following
observation estimate at one point in time:

\begin{prop}
\label{propsi} Under the assumptions of Theorem \ref{t11}, there exist $c>1$
and $M>1$ (both depending on $|x_{0}|$, $r$, $\left\Vert \left(  a_{0}%
,b_{0}\right)  \right\Vert _{L^{3}\left(  \Omega\right)  }$, $||k||_{L^{\infty
}\left(  \Omega\times\mathbb{R}_{+}\right)  }$, $d_{1}$, $d_{2}$) such that
for any $t>t_{1}\geq0$,
\[%
\begin{array}
[c]{ll}%
\left(  \left\Vert \left(  a-1,b-1\right)  \left(  \cdot,t\right)  \right\Vert
_{\left(  L^{2}\left(  \Omega\right)  \right)  ^{2}}^{2}\right)  ^{1+M} &
\leq\displaystyle e^{c\left(  1+\frac{1}{t-t_{1}}\right)  }\left\Vert \left(
a-1,b-1\right)  \left(  \cdot,t\right)  \right\Vert _{\left(  L^{2}\left(
B\left(  x_{0},r\right)  \right)  \right)  ^{2}}^{2}\\
& \quad\times\left(  \left\Vert \left(  a-1,b-1\right)  \left(  \cdot
,t_{1}\right)  \right\Vert _{\left(  L^{2}\left(  \Omega\right)  \right)
^{2}}^{2}\right)  ^{M}\text{ .}%
\end{array}
\]

\end{prop}

The proof of Proposition \ref{propsi} is given at the beginning of Section
\ref{sec3}. We first proceed with the \medskip

\textbf{{End of the Proof of Theorem \ref{t11}}}: First, we write down the
energy identity (\ref{ape2}) for the quantities $u_{1}:=a-1$, $u_{2}:=b-1$:
\begin{equation}
\frac{1}{2}\frac{d}{dt}\left\Vert \left(  u_{1},u_{2}\right)  \right\Vert
_{\left(  L^{2}\left(  \Omega\right)  \right)  ^{2}}^{2}+d_{1}\int_{\Omega
}\left\vert \nabla u_{1}\right\vert ^{2}+d_{2}\int_{\Omega}\left\vert \nabla
u_{2}\right\vert ^{2}+\int_{\Omega}k\left(  a+b\right)  \left\vert u_{2}%
-u_{1}\right\vert ^{2}=0\text{ ,} \label{2.1}%
\end{equation}
and the identity
\begin{equation}
\forall t\geq0\text{ ,}\qquad\int_{\Omega}[u_{1}(\cdot,t)+u_{2}(\cdot
,t)]=0\text{ ,} \label{uu2}%
\end{equation}
which is a direct consequence of (\ref{hypin2}), (\ref{ape1}). \medskip

Second, by Poincar\'{e}-Wirtinger's inequality (using identity (\ref{uu2}),
and denoting by $C_{p}$ the corresponding constant), the assumption $k\left(
\cdot,t\right)  \geq k_{0}>0$ on $B\left(  x_{0},r\right)  $, and remembering
(\ref{pmax}), we see that
\[%
\begin{array}
[c]{ll}%
2\left\Vert \left(  u_{1},u_{2}\right)  \right\Vert _{\left(  L^{2}\left(
B\left(  x_{0},r\right)  \right)  \right)  ^{2}}^{2} & =\left\Vert u_{1}%
+u_{2}\right\Vert _{L^{2}\left(  B\left(  x_{0},r\right)  \right)  }%
^{2}+\left\Vert u_{1}-u_{2}\right\Vert _{L^{2}\left(  B\left(  x_{0},r\right)
\right)  }^{2}\\
& \leq\left\Vert u_{1}+u_{2}\right\Vert _{L^{2}\left(  \Omega\right)  }%
^{2}+\displaystyle\frac{1}{2B_{0}k_{0}}\displaystyle\int_{\Omega}k\left(
a+b\right)  \left\vert u_{2}-u_{1}\right\vert ^{2}\\
& \leq C_{p}\left\Vert \nabla\left(  u_{1}+u_{2}\right)  \right\Vert
_{L^{2}\left(  \Omega\right)  }^{2}+\displaystyle\frac{1}{2B_{0}k_{0}%
}\displaystyle\int_{\Omega}k\left(  a+b\right)  \left\vert u_{2}%
-u_{1}\right\vert ^{2}\\
& \leq4\beta_{1}\left(  d_{1}\displaystyle\int_{\Omega}\left\vert \nabla
u_{1}\right\vert ^{2}+d_{2}\displaystyle\int_{\Omega}\left\vert \nabla
u_{2}\right\vert ^{2}+\displaystyle\int_{\Omega}k\left(  a+b\right)
\left\vert u_{2}-u_{1}\right\vert ^{2}\right)  \text{ ,}%
\end{array}
\]
with $\beta_{1}:=$max$\left(  \frac{C_{p}}{2d_{1}},\frac{C_{p}}{2d_{2}}%
,\frac{1}{8B_{0}k_{0}}\right)  $. \medskip

Combining the above estimate with (\ref{2.1}) and Proposition \ref{propsi}, we
deduce that
\[
\beta_{1}\frac{d}{dt}\left\Vert \left(  u_{1},u_{2}\right)  \left(
\cdot,t\right)  \right\Vert _{\left(  L^{2}\left(  \Omega\right)  \right)
^{2}}^{2}+\frac{1}{e^{c\left(  1+\frac{1}{t-t_{1}}\right)  }}\frac{\left(
\left\Vert \left(  u_{1},u_{2}\right)  \left(  \cdot,t\right)  \right\Vert
_{\left(  L^{2}\left(  \Omega\right)  \right)  ^{2}}^{2}\right)  ^{1+M}%
}{\left(  \left\Vert \left(  u_{1},u_{2}\right)  \left(  \cdot,t_{1}\right)
\right\Vert _{\left(  L^{2}\left(  \Omega\right)  \right)  ^{2}}^{2}\right)
^{M}}\leq0\text{ ,}%
\]
which can be rewritten
\begin{equation}
\label{np1}\frac{e^{-c\left(  1+\frac{1}{t-t_{1}}\right)  }M/\beta_{1}%
}{\left(  \left\Vert \left(  u_{1},u_{2}\right)  \left(  \cdot,t_{1}\right)
\right\Vert _{\left(  L^{2}\left(  \Omega\right)  \right)  ^{2}}^{2}\right)
^{M}}\leq\frac{d}{dt}\left[  \left(  \left\Vert \left(  u_{1},u_{2}\right)
\left(  \cdot,t\right)  \right\Vert _{\left(  L^{2}\left(  \Omega\right)
\right)  ^{2}}^{2}\right)  ^{-M}\right]  \text{ .}%
\end{equation}
Integrating (\ref{np1}) over $\left(  t_{1}+1,t_{2}\right)  $ with
$t_{2}>t_{1}+1 \ge1$, we obtain
\begin{equation}
\label{np2}\frac{e^{-2c}M\left(  t_{2}-t_{1}-1\right)  /\beta_{1}}{\left(
\left\Vert \left(  u_{1},u_{2}\right)  \left(  \cdot,t_{1}\right)  \right\Vert
_{\left(  L^{2}\left(  \Omega\right)  \right)  ^{2}}^{2}\right)  ^{M}}%
\leq\frac{1}{\left(  \left\Vert \left(  u_{1},u_{2}\right)  \left(
\cdot,t_{2}\right)  \right\Vert _{\left(  L^{2}\left(  \Omega\right)  \right)
^{2}}^{2}\right)  ^{M}}-\frac{1}{\left(  \left\Vert \left(  u_{1}%
,u_{2}\right)  \left(  \cdot,t_{1}+1\right)  \right\Vert _{\left(
L^{2}\left(  \Omega\right)  \right)  ^{2}}^{2}\right)  ^{M}}\text{ .}%
\end{equation}
But $\left\Vert \left(  u_{1},u_{2}\right)  \left(  \cdot,t_{1}+1\right)
\right\Vert _{\left(  L^{2}\left(  \Omega\right)  \right)  ^{2}}^{2}%
\leq\left\Vert \left(  u_{1},u_{2}\right)  \left(  \cdot,t_{1}\right)
\right\Vert _{\left(  L^{2}\left(  \Omega\right)  \right)  ^{2}}^{2}$ thanks
to (\ref{2.1}). Therefore,
\begin{equation}
\label{np3}\left\Vert \left(  u_{1},u_{2}\right)  \left(  \cdot,t_{2}\right)
\right\Vert _{\left(  L^{2}\left(  \Omega\right)  \right)  ^{2}}^{2}%
\leq\left(  \frac{1}{1+e^{-2c}M\left(  t_{2}-t_{1}-1\right)  /\beta_{1}%
}\right)  ^{1/M}\left\Vert \left(  u_{1},u_{2}\right)  \left(  \cdot
,t_{1}\right)  \right\Vert _{\left(  L^{2}\left(  \Omega\right)  \right)
^{2}}^{2}\text{ .}%
\end{equation}
Now, choose $t_{1}= 2m$ and $t_{2}=2\left(  m+1\right)  $, where
$m\in\mathbb{N}$ (so that $t_{2}>t_{1}+1 \ge1$). Then estimate (\ref{np3})
becomes
\begin{equation}
\label{np4}\left\Vert \left(  u_{1},u_{2}\right)  \left(  \cdot,2\left(
m+1\right)  \right)  \right\Vert _{\left(  L^{2}\left(  \Omega\right)
\right)  ^{2}}^{2}\leq\theta\left\Vert \left(  u_{1},u_{2}\right)  \left(
\cdot, 2m\right)  \right\Vert _{\left(  L^{2}\left(  \Omega\right)  \right)
^{2}}^{2}\text{ ,}%
\end{equation}
where $\theta:=\left(  \frac{1}{1+e^{-2c}M / \beta_{1}}\right)  ^{1/M}%
\in\left(  0,1\right)  $. A direct induction shows that
\begin{equation}
\label{np5}\left\Vert \left(  u_{1},u_{2}\right)  \left(  \cdot,2m\right)
\right\Vert _{\left(  L^{2}\left(  \Omega\right)  \right)  ^{2}}^{2}\leq
\theta^{m}\left\Vert \left(  u_{1},u_{2}\right)  \left(  \cdot,0\right)
\right\Vert _{\left(  L^{2}\left(  \Omega\right)  \right)  ^{2}}^{2}\text{ .}%
\end{equation}
Choosing $2m\leq t<2\left(  m+1\right)  $, we obtain thanks to (\ref{2.1})
that
\begin{equation}
\label{np6}\left\Vert \left(  u_{1},u_{2}\right)  \left(  \cdot,t\right)
\right\Vert _{\left(  L^{2}\left(  \Omega\right)  \right)  ^{2}}^{2}\leq
\frac{1}{\theta}e^{-t\frac{\left\vert \text{ln}\theta\right\vert }{2}%
}\left\Vert \left(  u_{1},u_{2}\right)  \left(  \cdot,0\right)  \right\Vert
_{\left(  L^{2}\left(  \Omega\right)  \right)  ^{2}}^{2}\text{ .}%
\end{equation}

We conclude the proof of Theorem \ref{t11} by taking $\gamma:=1/\theta$,
$\beta:=|\ln \theta|/2$.

\section{Heat systems with Neumann boundary conditions}

\label{sec3}

\subsection{An observation estimate}

\label{new31}

We consider in this section a system of two heat equations with homogeneous
Neumann boundary conditions and a second member (still on the ball $\Omega$
defined earlier):%
\begin{equation}
\left\{
\begin{array}
[c]{ll}%
{\partial}_{t}u_{1}-d_{1}\Delta u_{1}=v_{1}\text{ ,} & \quad\text{in}%
~\Omega\times\left(  0,+\infty\right)  \text{ ,}\\
{\partial}_{t}u_{2}-d_{2}\Delta u_{2}=v_{2}\text{ ,} & \quad\text{in}%
~\Omega\times\left(  0,+\infty\right)  \text{ ,}\\
\partial_{n}u_{1}=\partial_{n}u_{2}=0\text{ ,} & \quad\text{on}~\partial
\Omega\times\left(  0,+\infty\right)  \text{ .}%
\end{array}
\right.  \label{hs}%
\end{equation}
We will focus on functions $u_{1},u_{2}$, $v_{1},v_{2}$ satisfying the
following bounds, for some $K_{0}>0$,
\begin{equation}
\forall\left(  x,t\right)  \in\Omega\times\left(  0,+\infty\right)  \text{
,}\qquad\left\vert \left(  v_{1},v_{2}\right)  \left(  x,t\right)  \right\vert
^{2}\leq K_{0}\,\left(  \left\vert \left(  u_{1},u_{2}\right)  \left(
x,t\right)  \right\vert ^{2}+\left\vert \left(  u_{1},u_{2}\right)  \left(
x,t\right)  \right\vert ^{4}\right)  \text{ ,} \label{3.1}%
\end{equation}%
\begin{equation}
\forall i\in\left\{  1,2\right\}  \text{ , }\forall t\in\mathbb{R}_{+}\text{
,}\qquad\left\Vert u_{i}\left(  \cdot,t\right)  \right\Vert _{L^{3}\left(
\Omega\right)  }^{2}\leq K_{0}\text{ ,} \label{3.2}%
\end{equation}
and%
\begin{equation}
\forall\text{ }0\leq t_{1}<t_{2}\text{ , }\qquad\left\Vert \left(  u_{1}%
,u_{2}\right)  \left(  \cdot,t_{2}\right)  \right\Vert _{\left(  L^{2}\left(
\Omega\right)  \right)  ^{2}}^{2}\leq\left\Vert \left(  u_{1},u_{2}\right)
\left(  \cdot,t_{1}\right)  \right\Vert _{\left(  L^{2}\left(  \Omega\right)
\right)  ^{2}}^{2}\text{ .} \label{3.3}%
\end{equation}
\medskip

Our aim is to prove the following observation estimate at one point in time:

\begin{theo}
\label{mt2} Let $\Omega:=B\left(  0,R\right)  :=\left\{  x\in\mathbb{R}%
^{\normalfont{\text{n}}};\left\vert x\right\vert <R\right\}  $, when
$\normalfont{\text{n}}\in\{1,2,3\}$, be the centered ball of $\mathbb{R}%
^{\normalfont{\text{n}}}$ of measure $1$, and let $d_{1},d_{2}>0$. We also
consider $x_{0}\in\Omega$ and $r>0$, such that $B\left(  x_{0},r\right)
\subset\Omega$. Let $(u_{1},u_{2},v_{1},v_{2})$ be smooth ($C^{2}%
({\overline{\Omega}}\times\lbrack0,+\infty))$) solutions to system (\ref{hs})
satisfying estimates (\ref{3.1}), (\ref{3.2}) and (\ref{3.3}).

Then there exist $c>1$ and $M>1$ (both depending on $K_{0}$, $|x_{0}|$, $r$,
and $d_{1}$, $d_{2}$) such that for any $T>0$,%
\[%
\begin{array}
[c]{ll}%
\left(  \left\Vert \left(  u_{1},u_{2}\right)  \left(  \cdot,T\right)
\right\Vert _{\left(  L^{2}\left(  \Omega\right)  \right)  ^{2}}^{2}\right)
^{1+M} & \leq e^{c\left(  1+\frac{1}{T}\right)  }\left\Vert \left(
u_{1},u_{2}\right)  \left(  \cdot,T\right)  \right\Vert _{\left(  L^{2}\left(
B\left(  x_{0},r\right)  \right)  \right)  ^{2}}^{2}\\
& \quad\times\left(  \left\Vert \left(  u_{1},u_{2}\right)  \left(
\cdot,0\right)  \right\Vert _{\left(  L^{2}\left(  \Omega\right)  \right)
^{2}}^{2}\right)  ^{M}\text{ .}%
\end{array}
\]

\end{theo}

\bigskip

In order to make appear the subset $B\left(  x_{0},r\right)  $ in the right
hand side of the main estimate in Theorem 3.1, it is usual to use either
Carleman inequality or a logarithmic convexity method. We refer to \cite{Le
B}, \cite{HZ} for the strategy with Carleman inequalities. The observation
estimate in Theorem 3.1 is obtained in this work by logarithmic convexity (see
\cite{PW}, \cite{P} in the case of Dirichlet Boundary condition and \cite{BuP} for more
general geometries). Here we have chosen to make the study in the ball and
track the constants.  Originally, logarithmic convexity was used for getting backward estimates
(see \cite{BT}). Here, we have slightly deformed the process with a weight
function (see Proposition~\ref{fas}) and the gain in time in the left hand
side of estimate (\ref{f1}) allows to play and to make appear the subset
$B\left(  x_{0},r\right)  $.
 \medskip

The proof of Theorem
\ref{mt2} is divided into 5 steps, described in subsections \ref{new32} to
\ref{new36}.  We first show here that this Theorem enables to
write down the
\medskip

\textbf{{Proof of Proposition \ref{propsi}}}: Under the assumptions of
Proposition \ref{propsi} (which are those of Theorem \ref{t11}), we see that
$u_{1} := a-1$ and $u_{2} :=b-1$ are smooth ($C^{2}({\overline{\Omega}}
\times[0, +\infty))$) and satisfy the system
\[
\left\{
\begin{array}
[c]{ll}%
{\partial}_{t}u_{1}-d_{1}\Delta u_{1}=k\left(  u_{1} + u_{2} + 2\right)
\left(  u_{2}-u_{1}\right)  \text{ ,} & \quad\text{in}~\Omega\times\left(
0,+\infty\right)  \text{ ,}\\
{\partial}_{t}u_{2}-d_{2}\Delta u_{2}=-k\left(  u_{1} + u_{2} + 2\right)
\left(  u_{2}-u_{1}\right)  \text{ ,} & \quad\text{in}~\Omega\times\left(
0,+\infty\right)  \text{ ,}\\
\partial_{n}u_{1}=\partial_{n}u_{2}=0\text{ ,} & \quad\text{on}~\partial
\Omega\times\left(  0,+\infty\right)  \text{ .}%
\end{array}
\right.
\]
Considering $v_{1} := k\,(u_{1} + u_{2} + 2)\,(u_{2} - u_{1})$, $v_{2} :=
-k\,(u_{1} + u_{2} + 2)\,(u_{2} - u_{1})$, we see that $v_{1}$ and $v_{2}$ are
smooth ($C^{2}({\overline{\Omega}} \times[0, +\infty))$) and that (\ref{3.3})
holds thanks to (\ref{2.1}). \medskip

Moreover,
\[
|(v_{1},v_{2})|^{2}=2k^{2}\,(u_{1}+u_{2}+2)^{2}\,(u_{2}-u_{1})^{2}%
\leq32\left\Vert k\right\Vert _{L^{\infty}\left(  \Omega\times\mathbb{R}%
_{+}\right)  }^{2}(|(u_{1},u_{2})|^{2}+|(u_{1},u_{2})|^{4})\text{ ,}%
\]
and finally (using (\ref{fp1})),%
\[
\left\{
\begin{array}
[c]{ll}%
||(u_{1}(\cdot,t)||_{L^{3}\left(  \Omega\right)  }^{3}=\displaystyle\int
_{\Omega}\left\vert a(\cdot,t)-1\right\vert ^{3}\leq4\displaystyle\int
_{\Omega}a^{3}(\cdot,t)+4\leq4\displaystyle\int_{\Omega}\left(  a_{0}%
^{3}+b_{0}^{3}\right)  +4\text{ ,} & \\
||(u_{2}(\cdot,t)||_{L^{3}\left(  \Omega\right)  }^{3}=\displaystyle\int
_{\Omega}\left\vert b(\cdot,t)-1\right\vert ^{3}\leq4\displaystyle\int
_{\Omega}b^{3}(\cdot,t)+4\leq4\displaystyle\int_{\Omega}\left(  a_{0}%
^{3}+b_{0}^{3}\right)  +4\text{ ,} &
\end{array}
\right.
\]
so that (\ref{3.1}) and (\ref{3.2}) hold with $K_{0}:=$max$\left(  \left[
4\int_{\Omega}\left(  a_{0}^{3}+b_{0}^{3}\right)  +4\right]  ^{2/3}%
,32\left\Vert k\right\Vert _{L^{\infty}\left(  \Omega\times\mathbb{R}%
_{+}\right)  }^{2}\right)  $. \medskip

Then it is possible to use Theorem \ref{mt2} and to get Proposition
\ref{propsi} after a simple time translation. \medskip

Subsections \ref{new32} to \ref{new36} are devoted to the proof of Theorem
\ref{mt2}.

\subsection{Step 1: Change of functions}

\label{new32}

\bigskip

We recall that $\Omega$ is the centered ball of measure~$1$ of $\mathbb{R}%
^{\text{n}}$, with n$=1,2,3$. Without loss of generality, we suppose that
$x_{0}=(|x_{0}|,0,..,0)$. Then we consider
\begin{equation}
\psi\left(  x\right)  :=\psi(x_{1},..,x_{\text{n}})=\left(  R^{2}-\left\vert
x\right\vert ^{2}\right)  \,\left(  \frac{2|x_{0}|\,R}{|x_{0}|^{2}%
+R^{2}-2|x_{0}|x_{1}}\right)  \text{ ,}\label{defpsi}%
\end{equation}
which is well defined and is $C^{\infty}$ on an open ball containing
$\overline{\Omega}$. Moreover $\psi>0$ on $\Omega$, and $\psi=0$ on
$\partial\Omega$. One can check that
\[
\frac{\partial\psi}{\partial x_{k}}(x)=\frac{-4|x_{0}|R\,x_{k}}{|x_{0}%
|^{2}+R^{2}-2|x_{0}|x_{1}}\qquad\text{if}\qquad k\neq1\text{ ,}%
\]%
\[
\frac{\partial\psi}{\partial x_{1}}(x)=\frac{-4|x_{0}|R\,x_{1}}{|x_{0}%
|^{2}+R^{2}-2|x_{0}|x_{1}}+(R^{2}-|x|^{2})\frac{4|x_{0}|^{2}R}{(|x_{0}%
|^{2}+R^{2}-2|x_{0}|x_{1})^{2}}\text{ .}%
\]
Then it is easy to see that $\psi$ has a unique critical point at $x_{0}$ on
$\overline{\Omega}$, which is a global nondegenerate maximum. Indeed (for
$j,k=1,..,$n),
\[
\frac{\partial^{2}\psi}{\partial x_{j}\partial x_{k}}(x_{0})=-\frac{4|x_{0}%
|R}{R^{2}-|x_{0}|^{2}}\,\delta_{jk}\text{ .}%
\]

In particular $\partial_{n}\psi\leq0$ on $\partial\Omega$, and further, there
exist $c_{01},c_{02}>0$, depending only on $|x_{0}|$, such that for any $x$ in
a small neighborhood (also only depending on $|x_{0}|$) of $x_{0}$, the
following estimate holds:
\begin{equation}
c_{01}\left\vert \nabla\psi\left(  x\right)  \right\vert ^{2}\leq\psi\left(
x_{0}\right)  -\psi\left(  x\right)  \leq c_{02}\left\vert \nabla\psi\left(
x\right)  \right\vert ^{2}\text{ .}\label{gpp}%
\end{equation}
\bigskip

We introduce $f:=\left(  f_{i}\right)  _{1\leq i\leq4}$, where $f_{i}%
:=u_{i}\,e^{\Phi_{i}/2}$, and $\Phi_{i}\left(  x,t\right)  :=\frac
{s\varphi_{i}\left(  x\right)  }{\Gamma\left(  t\right)  }$, $s\in\left(
0,1\right]  $, $h\in\left(  0,1\right]  $,
\begin{equation}
\left\{
\begin{array}
[c]{ll}%
\Gamma\left(  t\right)  =T-t+h\text{ ,} & \quad\text{for any}~t\in\left[
0,T\right]  \text{ ,}\\
\varphi_{1}\left(  x\right)  =\psi\left(  x\right)  -\psi\left(  x_{0}\right)
\text{ ,} & \quad\text{for any}~x\in\overline{\Omega}\text{ ,}\\
\varphi_{3}\left(  x\right)  =-\psi\left(  x\right)  -\psi\left(
x_{0}\right)  \text{ ,} & \quad\text{for any}~x\in\overline{\Omega}\text{ ,}%
\end{array}
\right.  \label{gaphi}%
\end{equation}
and
\begin{equation}
\varphi_{2}:=\varphi_{1}\leq0\text{ ,}\quad\varphi_{4}:=\varphi_{3}~(\text{so
that{ }}\Phi_{2}:=\Phi_{1}\leq0\text{ ,}\quad\Phi_{4}:=\Phi_{3})\text{,}\quad
u_{3}:=u_{1}\text{ ,}\quad u_{4}:=u_{2}\text{ .}\label{prodo0}%
\end{equation}
\medskip

Notice that on $\partial\Omega\times(0,T)$,
\begin{equation}
\varphi_{1}=\varphi_{3}\text{ ,}\qquad\partial_{t}\varphi_{1}=\partial
_{t}\varphi_{3}\text{ ,}\qquad\partial_{n}\varphi_{1}+\partial_{n}\varphi
_{3}=0\text{ ,}\label{prodo1}%
\end{equation}%
\begin{equation}
\varphi_{2}=\varphi_{4}\text{ ,}\qquad\partial_{n}\varphi_{2}+\partial
_{n}\varphi_{4}=0\text{ ,}\label{prodo2}%
\end{equation}
so that, still on $\partial\Omega\times(0,T)$,
\begin{equation}
\Phi_{1}=\Phi_{3}\text{ ,}\qquad\partial_{t}\Phi_{1}=\partial_{t}\Phi
_{3}\text{ ,}\qquad\partial_{n}\Phi_{1}+\partial_{n}\Phi_{3}=0\text{
,}\label{prodo3}%
\end{equation}%
\begin{equation}
\Phi_{2}=\Phi_{4}\text{ ,}\qquad\partial_{n}\Phi_{2}+\partial_{n}\Phi
_{4}=0\text{ .}\label{prodo4}%
\end{equation}

We look for the equation solved by $f_{i}$ by computing $e^{\Phi_{i}/2}\left(
\partial_{t}-d_{i}\Delta\right)  \left(  e^{-\Phi_{i}/2}f_{i}\right)  $,
where
\begin{equation}
d_{3}:=d_{1}\qquad\text{and}\qquad d_{4}:=d_{2}\text{ .}\label{dd1}%
\end{equation}
\medskip

We introduce for that purpose the operators
\begin{equation}
\label{as}\left\{
\begin{array}
[c]{ll}%
\mathcal{A}_{i}f_{i} :=-d_{i}\nabla\Phi_{i}\cdot\nabla f_{i}-\frac{1}{2}%
d_{i}\Delta\Phi_{i}f_{i}\text{ ,} & \\
\mathcal{S}_{i}f_{i} :=-d_{i}\Delta f_{i}-\eta_{i}f_{i} \text{ ,} &
\end{array}
\right.
\end{equation}
where $i=1,..,4$ and
\begin{equation}
\label{eta}\eta_{i} := \frac{1}{2}\partial_{t}\Phi_{i}+\frac{1}{4}%
d_{i}\left\vert \nabla\Phi_{i}\right\vert ^{2}\text{ .}%
\end{equation}

We also define $\mathcal{S}f:=\left(  \mathcal{S}_{i}f_{i}\right)  _{1\leq
i\leq4}$, $\mathcal{A}f:=\left(  \mathcal{A}_{i}f_{i}\right)  _{1\leq i\leq4}%
$, and $\digamma:=\left(  v_{i}e^{\Phi_{i}/2}\right)  _{1\leq i\leq4}$ where
\begin{equation}
v_{3}:=v_{1}\text{ ,}\qquad v_{4}:=v_{2}\text{ .}\label{pv3}%
\end{equation}

After this change of functions and the introduction of the new notations, the
system (\ref{hs}) rewrites
\begin{equation}
\left\{
\begin{array}
[c]{ll}%
\partial_{t}f+\mathcal{S}f=\mathcal{A}f+\digamma\text{ ,} & \\
\partial_{n}f_{i}-\frac{1}{2}\partial_{n}\Phi_{i}f_{i}=0\text{ on }%
\partial\Omega\times\left(  0,T\right)  \text{ ,}\qquad i=1,..,4\text{ .} &
\end{array}
\right.  \label{nhs}%
\end{equation}

Let now $\left\langle \cdot,\cdot\right\rangle $ denote the usual scalar
product in $\left(  L^{2}\left(  \Omega\right)  \right)  ^{4}$, and
$\left\Vert \cdot\right\Vert $ be its corresponding norm. We regroup some
useful identities in the following:

\begin{prop}
\label{fas} For any (smooth enough, $\mathbb{R}^{4}$-valued) functions
$f:=f(x,t)=(f_{i})_{i=1,..,4}$, any constants $d_{i}>0$, and any (smooth
enough) functions $\Phi_{i}:=\Phi_{i}(x,t)$, $u_{i}=f_{i}\,e^{-\Phi_{i}/2}$
($i=1,..,4$), the following identities hold as soon as (\ref{prodo0}) --
(\ref{nhs}) hold:
\begin{equation}
\left\{
\begin{array}
[c]{ll}%
\left\langle \mathcal{A}f,f\right\rangle =0\text{ ,} & \\
\left\langle \mathcal{S}f,f\right\rangle =\displaystyle\sum_{i=1,..,4}%
\bigg[d_{i}\displaystyle\int_{\Omega}\left\vert \nabla f_{i}\right\vert
^{2}-\displaystyle\int_{\Omega}\eta_{i}\left\vert f_{i}\right\vert
^{2}\bigg]\text{ ,} & \\
\displaystyle\frac{d}{dt}\left\langle \mathcal{S}f,f\right\rangle
=\displaystyle\sum_{i=1,..,4}\displaystyle\int_{\Omega}\left(  -\partial
_{t}\eta_{i}\right)  \left\vert f_{i}\right\vert ^{2}+2\left\langle
\mathcal{S}f,\partial_{t}f\right\rangle :=\left\langle \mathcal{S}^{\prime
}f,f\right\rangle +2\left\langle \mathcal{S}f,\mathcal{\partial}%
_{t}f\right\rangle \text{ .} &
\end{array}
\right.  \label{asds}%
\end{equation}

\end{prop}

\bigskip

Proof of Proposition~\ref{fas} .- Thanks to an integration by parts,
\[%
\begin{array}
[c]{ll}
& \quad\left\langle \mathcal{A}f,f\right\rangle :=\displaystyle\sum
_{i=1,..,4}\displaystyle\int_{\Omega}\left(  -d_{i}\nabla\Phi_{i}\cdot\nabla
f_{i}-\frac{1}{2}d_{i}\Delta\Phi_{i}f_{i}\right)  f_{i}\\
& =-\displaystyle\sum_{i=1,..,4}\displaystyle\int_{\Omega} \bigg[ d_{i}%
\nabla\Phi_{i}\cdot\nabla\left(  \frac{1}{2}\left\vert f_{i}\right\vert
^{2}\right)  +\frac{1}{2}d_{i}\Delta\Phi_{i}\left\vert f_{i}\right\vert ^{2}
\bigg]\\
& =-\displaystyle\int_{\partial\Omega} \bigg[ d_{1}\partial_{n}\Phi_{1}\left(
\frac{1}{2}\left\vert f_{1}\right\vert ^{2}\right)  +d_{2}\partial_{n}\Phi
_{2}\left(  \frac{1}{2}\left\vert f_{2}\right\vert ^{2}\right)  +d_{3}%
\partial_{n}\Phi_{3}\left(  \frac{1}{2}\left\vert f_{3}\right\vert
^{2}\right)  +d_{4}\partial_{n}\Phi_{4}\left(  \frac{1}{2}\left\vert
f_{4}\right\vert ^{2}\right)  \bigg]\\
& =-\displaystyle\int_{\partial\Omega}d_{1}\left(  \partial_{n}\Phi_{1}%
e^{\Phi_{1}}+\partial_{n}\Phi_{3}e^{\Phi_{3}}\right)  \left(  \frac{1}%
{2}\left\vert u_{1}\right\vert ^{2}\right)  -\displaystyle\int_{\partial
\Omega}d_{2}\left(  \partial_{n}\Phi_{1}e^{\Phi_{1}}+\partial_{n}\Phi
_{3}e^{\Phi_{3}}\right)  \left(  \frac{1}{2}\left\vert u_{2}\right\vert
^{2}\right)
\end{array}
\]
using (\ref{prodo0}), (\ref{dd1}) and recalling that $f_{i}=u_{i}e^{\Phi
_{i}/2}$.

Now, (\ref{prodo3}) implies that
\begin{equation}
\label{prodo9}\partial_{n}\Phi_{1}e^{\Phi_{1}}+\partial_{n}\Phi_{3}e^{\Phi
_{3}}=0\text{ on }\partial\Omega\times\left(  0,T\right)  \text{ .}%
\end{equation}
This completes the proof of the identity $\left\langle \mathcal{A}%
f,f\right\rangle =0$. \medskip

Then, we observe that
\[%
\begin{array}
[c]{ll}
& \quad\left\langle \mathcal{S}f,f\right\rangle -\left[  \displaystyle\sum
_{i=1,..,4}d_{i}\displaystyle\int_{\Omega}\left\vert \nabla f_{i}\right\vert
^{2}-\displaystyle\sum_{i=1,..,4}\displaystyle\int_{\Omega}\eta_{i}\left\vert
f_{i}\right\vert ^{2}\right] \\
& =-\displaystyle\sum_{i=1,..,4}\displaystyle\int_{\partial\Omega}d_{i}%
f_{i}\partial_{n}f_{i}\\
& =-\displaystyle\sum_{i=1,..,4}\displaystyle\int_{\partial\Omega}d_{i}\left(
\frac{1}{2}\partial_{n}\Phi_{i}\left\vert f_{i}\right\vert ^{2}\right)  \text{
because of (\ref{nhs} }ii)\text{ ) }\\
& =0\text{ similarly as for }\left\langle \mathcal{A}f,f\right\rangle \text{
.}%
\end{array}
\]
We finally compute $\displaystyle\frac{d}{dt}\left\langle \mathcal{S}%
f,f\right\rangle :=\displaystyle\frac{d}{dt}\left(  \displaystyle\sum
_{i=1,..,4}\bigg[d_{i}\displaystyle\int_{\Omega}\left\vert \nabla
f_{i}\right\vert ^{2}-\displaystyle\int_{\Omega}\eta_{i}\left\vert
f_{i}\right\vert ^{2}\bigg]\right)  $. By an integration by parts,
\[%
\begin{array}
[c]{ll}%
\displaystyle\frac{d}{dt}\left\langle \mathcal{S}f,f\right\rangle  &
=\displaystyle\sum_{i=1,..,4}\bigg[d_{i}\displaystyle\int_{\Omega}2\nabla
f_{i}\cdot\nabla\partial_{t}f_{i}-\displaystyle\int_{\Omega}\partial_{t}%
\eta_{i}\left\vert f_{i}\right\vert ^{2}-\displaystyle\int_{\Omega}2\eta
_{i}f_{i}\partial_{t}f_{i}\bigg]\\
& =-\displaystyle\sum_{i=1,..,4}d_{i}\displaystyle\int_{\Omega}2\Delta
f_{i}\partial_{t}f_{i}+\displaystyle\sum_{i=1,..,4}\bigg[d_{i}%
\displaystyle\int_{\partial\Omega}2\partial_{n}f_{i}\partial_{t}%
f_{i}-\displaystyle\int_{\Omega}\partial_{t}\eta_{i}\left\vert f_{i}%
\right\vert ^{2}-\displaystyle\int_{\Omega}2\eta_{i}f_{i}\partial_{t}%
f_{i}\bigg]\\
& =\displaystyle\sum_{i=1,..,4}\displaystyle2\int_{\partial\Omega}%
d_{i}\partial_{n}f_{i}\partial_{t}f_{i}+\displaystyle\sum_{i=1,..,4}%
\displaystyle\int_{\Omega}\left(  -\partial_{t}\eta_{i}\right)  \left\vert
f_{i}\right\vert ^{2}+2\left\langle \mathcal{S}f,\partial_{t}f\right\rangle
\text{ .}%
\end{array}
\]
But
\[%
\begin{array}
[c]{ll}%
\displaystyle\sum_{i=1,..,4}\displaystyle\int_{\partial\Omega}d_{i}%
\partial_{n}f_{i}\partial_{t}f_{i} & =\displaystyle\sum_{i=1,..,4}%
\displaystyle\int_{\partial\Omega}d_{i}\left(  \frac{1}{2}\partial_{n}\Phi
_{i}f_{i}\partial_{t}f_{i}\right)  \text{ because of (\ref{nhs} }ii)\text{
)}\\
& =\displaystyle\sum_{i=1,..,4}\displaystyle\frac{1}{2}\int_{\partial\Omega
}d_{i}\partial_{n}\Phi_{i}u_{i}e^{\Phi_{i}/2}\left(  \partial_{t}u_{i}%
e^{\Phi_{i}/2}+u_{i}\frac{1}{2}\partial_{t}\Phi_{i}e^{\Phi_{i}/2}\right)
\text{ .}%
\end{array}
\]
The first contribution is proportional to
\[%
\begin{array}
[c]{ll}
& \quad\displaystyle\int_{\partial\Omega}\bigg[d_{1}\partial_{n}\Phi_{1}%
u_{1}\partial_{t}u_{1}e^{\Phi_{1}}+d_{2}\partial_{n}\Phi_{2}u_{2}\partial
_{t}u_{2}e^{\Phi_{2}}+d_{3}\partial_{n}\Phi_{3}u_{3}\partial_{t}u_{3}%
e^{\Phi_{3}}+d_{4}\partial_{n}\Phi_{4}u_{4}\partial_{t}u_{4}e^{\Phi_{4}%
}\bigg]\\
& =\displaystyle\int_{\partial\Omega}\bigg[d_{1}\left(  \partial_{n}\Phi
_{1}e^{\Phi_{1}}+\partial_{n}\Phi_{3}e^{\Phi_{3}}\right)  u_{1}\partial
_{t}u_{1}+d_{2}\left(  \partial_{n}\Phi_{1}e^{\Phi_{1}}+\partial_{n}\Phi
_{3}e^{\Phi_{3}}\right)  u_{2}\partial_{t}u_{2}\bigg]
\end{array}
\]
using (\ref{prodo0}), (\ref{dd1}). Thanks to (\ref{prodo9}), we see that
$\displaystyle\sum_{i=1,..,4}\displaystyle\int_{\partial\Omega}d_{i}%
\partial_{n}\Phi_{i}u_{i}\partial_{t}u_{i}e^{\Phi_{i}}=0$.

The last contribution is proportional to
\[%
\begin{array}
[c]{ll}
& \quad\displaystyle\sum_{i=1,..,4}\displaystyle\int_{\partial\Omega}%
d_{i}\partial_{n}\Phi_{i}\partial_{t}\Phi_{i}\left\vert u_{i}\right\vert
^{2}e^{\Phi_{i}}\\
& =\displaystyle\int_{\partial\Omega}\bigg[d_{1}\left(  \partial_{n}\Phi
_{1}\partial_{t}\Phi_{1}e^{\Phi_{1}}+\partial_{n}\Phi_{3}\partial_{t}\Phi
_{3}e^{\Phi_{3}}\right)  \left\vert u_{1}\right\vert ^{2}+d_{2}\left(
\partial_{n}\Phi_{1}\partial_{t}\Phi_{1}e^{\Phi_{1}}+\partial_{n}\Phi
_{3}\partial_{t}\Phi_{3}e^{\Phi_{3}}\right)  \left\vert u_{2}\right\vert
^{2}\bigg]\\
& =0
\end{array}
\]
where in the second line, we used (\ref{prodo0}), (\ref{dd1}). In the third
line, we used the identity $\partial_{t}\Phi_{1}=\partial_{t}\Phi_{3}$ on
$\partial\Omega\times\left(  0,T\right)  $, which is a consequence of
(\ref{prodo3}). \medskip

This completes the proof of Proposition \ref{fas}.

\subsection{Step 2: Energy estimates}

Thanks to a standard energy method, we get, starting from eq. (\ref{nhs}) and
using the first identity of Proposition \ref{fas}:
\begin{equation}
\frac{1}{2}\frac{d}{dt}\left\Vert f\right\Vert ^{2}+\left\langle
\mathcal{S}f,f\right\rangle =\left\langle \digamma,f\right\rangle \text{ .}
\label{sem}%
\end{equation}
Introducing the frequency function
\begin{equation}
N:=N\left(  t\right)  =\frac{\left\langle \mathcal{S}f,f\right\rangle
}{\left\Vert f\right\Vert ^{2}}\text{ ,} \label{freq}%
\end{equation}
we see that, thanks to (\ref{nhs}), to the third part of Proposition \ref{fas},
and to identity (\ref{sem}),
\[%
\begin{array}
[c]{ll}%
N^{\prime}\left(  t\right)  \left(  \left\Vert f\right\Vert ^{2}\right)  ^{2}
& =\left(  \left\langle \mathcal{S}^{\prime}f,f\right\rangle +2\left\langle
\mathcal{S}f,\partial_{t}f\right\rangle \right)  \left\Vert f\right\Vert
^{2}-\left\langle \mathcal{S}f,f\right\rangle \left(  -2\left\langle
\mathcal{S}f,f\right\rangle +2\left\langle \digamma,f\right\rangle \right) \\
& =\left(  \left\langle \mathcal{S}^{\prime}f,f\right\rangle +2\left\langle
\mathcal{S}f,\mathcal{A}f\right\rangle \right)  \left\Vert f\right\Vert
^{2}-2\left\Vert \mathcal{S}f\right\Vert ^{2}\left\Vert f\right\Vert
^{2}+2\left\langle \mathcal{S}f,\digamma\right\rangle \left\Vert f\right\Vert
^{2}\\
& \quad+2\left\langle \mathcal{S}f,f\right\rangle ^{2}-2\left\langle
\mathcal{S}f,f\right\rangle \left\langle \digamma,f\right\rangle \\
& =\left(  \left\langle \mathcal{S}^{\prime}f,f\right\rangle +2\left\langle
\mathcal{S}f,\mathcal{A}f\right\rangle \right)  \left\Vert f\right\Vert
^{2}-2\left\Vert \mathcal{S}f-\frac{1}{2}\digamma\right\Vert ^{2}\left\Vert
f\right\Vert ^{2}+\frac{1}{2}\left\Vert \digamma\right\Vert ^{2}\left\Vert
f\right\Vert ^{2}\\
& \quad+2\left\langle \mathcal{S}f-\frac{1}{2}\digamma,f\right\rangle
^{2}-\frac{1}{2}\left\langle \digamma,f\right\rangle ^{2}\text{ .}%
\end{array}
\]
Thanks to Cauchy-Schwarz inequality, we obtain the following estimate for
$N^{\prime}\left(  t\right)  $:
\begin{equation}
N^{\prime}\left(  t\right)  \left\Vert f\right\Vert ^{2}\leq\left\langle
\mathcal{S}^{\prime}f,f\right\rangle +2\left\langle \mathcal{S}f,\mathcal{A}%
f\right\rangle +\left\Vert \digamma\right\Vert ^{2}\text{ .} \label{ineqnpr}%
\end{equation}

\subsection{Step 3: Carleman commutator estimates}

The next ingredient in the proof of Theorem \ref{mt2} is the following:

\begin{prop}
\label{prop2} Under the assumptions of Theorem \ref{mt2} (and using the notations
(\ref{gaphi}) -- (\ref{eta})), there exist $s_{2}\in(0,1]$, $C_{0}\in\left(
0,1\right)  $, $C_{1}>1$ depending only on $K_{0}$, $|x_{0}|$, $d_{1}$,
$d_{2}$, such that when $s\in\left(  0,s_{2}\right]  $, $h\in(0,1]$,

\begin{enumerate}
\item[$\left(  i\right)  $] \qquad\qquad$\eta_{i}\leq0$ \, and \,
$\left\langle \mathcal{S}f,f\right\rangle \geq0$,

\item[$\left(  ii\right)  $]
\begin{equation}
\left\Vert \digamma\right\Vert ^{2}\leq C_{1}\left(  \left\Vert f\right\Vert
^{2}+\left\langle \mathcal{S}f,f\right\rangle \right)  \text{ ,} \label{p22}%
\end{equation}

\item[$\left(  iii\right)  $]
\begin{equation}
\left\langle \mathcal{S}^{\prime}f,f\right\rangle +2\left\langle
\mathcal{S}f,\mathcal{A}f\right\rangle \leq\frac{1+C_{0}}{\Gamma}\left\langle
\mathcal{S}f,f\right\rangle +\frac{C_{1}}{h^{2}}\left\Vert f\right\Vert
^{2}\text{ .} \label{p23}%
\end{equation}

\end{enumerate}
\end{prop}

The proof of Proposition \ref{prop2} is postponed to Section \ref{sec4}.

\subsection{Step 4: Use of a differential inequality}

We first observe that thanks to (\ref{sem}), Young's inequality and
(\ref{p22}), when $h \in(0,1]$, we get (for $s \in(0, s_{2}]$)%

\begin{equation}
\left\vert \frac{1}{2}\frac{d}{dt}\left\Vert f\right\Vert ^{2}+\left\langle
\mathcal{S}f,f\right\rangle \right\vert =\left\vert \left\langle
\digamma,f\right\rangle \right\vert \leq\frac{1}{2C_{1}}\left\Vert
\digamma\right\Vert ^{2}+\frac{C_{1}}{2}\left\Vert f\right\Vert ^{2}\leq
\frac{1}{2}\left\langle \mathcal{S}f,f\right\rangle +\frac{C_{1}}{h}\left\Vert
f\right\Vert ^{2}\text{ .} \label{hy21}%
\end{equation}
Moreover, thanks to (\ref{ineqnpr}), (\ref{p22}) and (\ref{p23}), we also get
(when $h\in(0,1]$):
\begin{equation}
N^{\prime}\leq\frac{\left\langle \mathcal{S}^{\prime}f,f\right\rangle
+2\left\langle \mathcal{S}f,\mathcal{A}f\right\rangle +\left\Vert
\digamma\right\Vert ^{2}}{\left\Vert f\right\Vert ^{2}}\leq\left(
\frac{1+C_{0}}{\Gamma}+C_{1}\right)  N+\frac{2C_{1}}{h^{2}}\text{ .}
\label{hy22}%
\end{equation}

\bigskip

We now wish to use the following Lemma, proven in \cite{BP}:

\bigskip

\begin{lem}
\label{lembp} Let $h>0$, $T>0$\ and $F_{1},F_{2}\in C\left(  \left[
0,T\right]  \right)  $. Consider two nonnegative functions $y,N\in
C^{1}\left(  \left[  0,T\right]  \right)  $\ such that (when $t\in\lbrack
0,T)$):
\[
\left\{
\begin{array}
[c]{ll}%
\left\vert \displaystyle\frac{1}{2}y^{\prime}\left(  t\right)  +N\left(
t\right)  y\left(  t\right)  \right\vert \leq\left(  \displaystyle\frac{1}%
{2}N\left(  t\right)  +\displaystyle\frac{C_{0}}{T-t+h}+C_{1}\right)  y\left(
t\right)  +F_{1}\left(  t\right)  y\left(  t\right)  \text{ ,} & \\
N^{\prime}\left(  t\right)  \leq\left(  \displaystyle\frac{1+C_{0}}%
{T-t+h}+C_{1}\right)  N\left(  t\right)  +F_{2}\left(  t\right)  \text{ ,} &
\end{array}
\right.
\]
where $C_{0},C_{1}\geq0$.

Then for any $0\leq t_{1}<t_{2}<t_{3}\leq T$, the following estimate holds:
\begin{equation}
y\left(  t_{2}\right)  ^{1+M}\leq e^{D}\left(  \frac{T-t_{1}+h}{T-t_{3}%
+h}\right)  ^{3C_{0}\left(  1+M\right)  }\,y\left(  t_{3}\right)  y\left(
t_{1}\right)  ^{M}\text{ ,} \label{y3}%
\end{equation}
with
\[
M:=3\frac{\displaystyle\int_{t_{2}}^{t_{3}}\frac{e^{tC_{1}}}{\left(
T-t+h\right)  ^{1+C_{0}}}dt}{\displaystyle\int_{t_{1}}^{t_{2}}\frac{e^{tC_{1}%
}}{\left(  T-t+h\right)  ^{1+C_{0}}}dt}%
\]
and
\[
D:=3\left(  1+M\right)  \left[  \left(  t_{3}-t_{1}\right)  \left(  C_{1}%
+\int_{t_{1}}^{t_{3}}\left\vert F_{2}\right\vert dt\right)  +\int_{t_{1}%
}^{t_{3}}\left\vert F_{1}\right\vert dt\right]  \text{ .}%
\]

\end{lem}

Consider now $h\in\left(  0,1\right]  $\ and $\ell>1$\ such that $\ell
h<$min$\left(  1/2,T/4\right)  $. Taking $t_{3}:=T$, $t_{2}:=T-\ell h$, and
$t_{1}:=T-2\ell h$, in Lemma \ref{lembp}, estimate (\ref{y3}) becomes
\[
y\left(  T-\ell h\right)  ^{1+M_{\ell}}\leq e^{D_{\ell}}\left(  2\ell
+1\right)  ^{3C_{0}\left(  1+M_{\ell}\right)  }\,y\left(  T\right)  y\left(
T-2\ell h\right)  ^{M_{\ell}}\text{ ,}%
\]
where
\begin{equation}
D_{\ell}:=3\left(  1+M_{\ell}\right)  \left(  C_{1}+\displaystyle\int_{T-2\ell
h}^{T}\left(  \left\vert F_{1}\right\vert +2\ell h\left\vert F_{2}\right\vert
\right)  dt\right)  \text{ ,} \label{defdl}%
\end{equation}
and
\begin{equation}
M_{\ell}:=3\frac{\displaystyle\int_{T-\ell h}^{T}\frac{e^{tC_{1}}}{\left(
T-t+h\right)  ^{1+C_{0}}}dt}{\displaystyle\int_{T-2\ell h}^{T-\ell h}%
\frac{e^{tC_{1}}}{\left(  T-t+h\right)  ^{1+C_{0}}}dt}\,\text{.} \label{defml}%
\end{equation}

\bigskip

We now observe that thanks to inequalities (\ref{hy21}), (\ref{hy22}), the
assumptions of Lemma~\ref{lembp} are fulfilled when $y\left(  t\right)
:=\left\Vert f\left(  \cdot,t\right)  \right\Vert ^{2}$, $N$ is the frequency
function given by (\ref{freq}), $h\in(0,1]$, $F_{1}(t):=\displaystyle\frac
{C_{1}}{h}$, $F_{2}(t):=\displaystyle\frac{2C_{1}}{h^{2}}$ and $C_{0}$,
$C_{1}$ are the constants appearing in Proposition \ref{prop2}. \medskip

Therefore, thanks to Lemma \ref{lembp}, for $h\in(0,1]$ and $\ell>1$ such that
$\ell h<$min$(1/2,T/4)$, the following estimate holds (for $s\in(0,s_{2}]$):
\begin{equation}
\left(  \left\Vert f\left(  \cdot,T-\ell h\right)  \right\Vert ^{2}\right)
^{1+M_{\ell}}\leq K_{\ell}\,\left(  \left\Vert f\left(  \cdot,T\right)
\right\Vert ^{2}\right)  \left(  \left\Vert f\left(  \cdot,T-2\ell h\right)
\right\Vert ^{2}\right)  ^{M_{\ell}}\text{ ,} \label{estiqua}%
\end{equation}
where $K_{\ell}:=e^{D_{\ell}}\left(  2\ell+1\right)  ^{3C_{0}\left(
1+M_{\ell}\right)  }$, with $D_{\ell}=3C_{1}\,\left(  1+M_{\ell}\right)
\left(  1+2\ell+ 8\ell^{2} \right)  $, and $M_{\ell}$ given by (\ref{defml}).

\subsection{Step 5: Introducing $B\left(  x_{0},r\right)  $}

\label{new36}

Observe that
\[
\left\Vert \left(  f_{1},f_{2}\right)  \right\Vert _{\left(  L^{2}\left(
\Omega\right)  \right)  ^{2}}^{2}\leq\left\Vert f\right\Vert ^{2}%
\leq2\left\Vert \left(  f_{1},f_{2}\right)  \right\Vert _{\left(  L^{2}\left(
\Omega\right)  \right)  ^{2}}^{2}\text{ ,}%
\]
since $\varphi_{3}-\varphi_{1}=-2\psi\leq0$ on $\Omega$. Therefore, thanks to
(\ref{estiqua}),
\begin{equation}%
\begin{array}
[c]{ll}
& \quad\left(  \left\Vert \left(  f_{1},f_{2}\right)  \left(  \cdot,T-\ell
h\right)  \right\Vert _{\left(  L^{2}\left(  \Omega\right)  \right)  ^{2}}%
^{2}\right)  ^{1+M_{\ell}}\\
& \leq K_{\ell}\,\left(  2\left\Vert \left(  f_{1},f_{2}\right)  \left(
\cdot,T\right)  \right\Vert _{\left(  L^{2}\left(  \Omega\right)  \right)
^{2}}^{2}\right)  \left(  2\left\Vert \left(  f_{1},f_{2}\right)  \left(
\cdot,T-2\ell h\right)  \right\Vert _{\left(  L^{2}\left(  \Omega\right)
\right)  ^{2}}^{2}\right)  ^{M_{\ell}}\text{ .}%
\end{array}
\label{f1}%
\end{equation}
First, notice that thanks to estimates (\ref{3.3}) and (\ref{prodo0}),
\begin{equation}
\left\Vert \left(  f_{1},f_{2}\right)  \left(  \cdot,T-2\ell h\right)
\right\Vert _{\left(  L^{2}\left(  \Omega\right)  \right)  ^{2}}^{2}%
\leq\left\Vert \left(  u_{1},u_{2}\right)  \left(  \cdot,0\right)  \right\Vert
_{\left(  L^{2}\left(  \Omega\right)  \right)  ^{2}}^{2}\text{ .} \label{f2}%
\end{equation}
Second, we make $B\left(  x_{0},r\right)  $ appear out of $\left\Vert \left(
f_{1},f_{2}\right)  \left(  \cdot,T\right)  \right\Vert _{\left(  L^{2}\left(
\Omega\right)  \right)  ^{2}}^{2}$ as follows:%
\begin{equation}%
\begin{array}
[c]{ll}%
\left\Vert \left(  f_{1},f_{2}\right)  \left(  \cdot,T\right)  \right\Vert
_{\left(  L^{2}\left(  \Omega\right)  \right)  ^{2}}^{2} & =\displaystyle\int
_{B\left(  x_{0},r\right)  }\left\vert \left(  u_{1},u_{2}\right)  \left(
\cdot,T\right)  \right\vert ^{2}e^{\frac{s}{h}\varphi_{1}}+\displaystyle\int
_{\Omega\left\backslash B\left(  x_{0},r\right)  \right.  }\left\vert \left(
u_{1},u_{2}\right)  \left(  \cdot,T\right)  \right\vert ^{2}e^{\frac{s}%
{h}\varphi_{1}}\\
& \leq\left\Vert \left(  u_{1},u_{2}\right)  \left(  \cdot,T\right)
\right\Vert _{\left(  L^{2}\left(  B\left(  x_{0},r\right)  \right)  \right)
^{2}}^{2}+e^{-\frac{s\mu_{0}}{h}}\left\Vert \left(  u_{1},u_{2}\right)
\left(  \cdot,0\right)  \right\Vert _{\left(  L^{2}\left(  \Omega\right)
\right)  ^{2}}^{2}%
\end{array}
\label{f3}%
\end{equation}
thanks to estimates (\ref{3.3}), (\ref{prodo0}), and the fact that on
$\Omega\left\backslash B\left(  x_{0},r\right)  \right.  $, $\varphi_{1}%
\leq-\mu_{0}$ for some $\mu_{0}>0$ depending only on $|x_{0}|$ and $r$.
\medskip

Third, thanks to estimate (\ref{3.3}),%
\begin{equation}%
\begin{array}
[c]{ll}%
\left\Vert \left(  u_{1},u_{2}\right)  \left(  \cdot,T\right)  \right\Vert
_{\left(  L^{2}\left(  \Omega\right)  \right)  ^{2}}^{2} & \leq
\displaystyle\int_{\Omega}\left\vert \left(  u_{1},u_{2}\right)  \left(
\cdot,T-\ell h\right)  \right\vert ^{2}e^{\frac{s}{\left(  \ell+1\right)
h}\varphi_{1}}e^{-\frac{s}{\left(  \ell+1\right)  h}\varphi_{1}}\\
& \leq e^{\frac{s\mu_{1}}{\left(  \ell+1\right)  h}}\left\Vert \left(
f_{1},f_{2}\right)  \left(  \cdot,T-\ell h\right)  \right\Vert _{\left(
L^{2}\left(  \Omega\right)  \right)  ^{2}}^{2}\text{ ,}%
\end{array}
\label{f4}%
\end{equation}
where $\mu_{1}:=\underset{x\in\overline{\Omega}}{\text{sup}}(-\varphi_{1}(x))$
only depends on $|x_{0}|$. \medskip

Combining estimates (\ref{f1}) to (\ref{f4}), we can deduce%
\begin{equation}%
\begin{array}
[c]{ll}
& \quad\left(  \left\Vert \left(  u_{1},u_{2}\right)  \left(  \cdot,T\right)
\right\Vert _{\left(  L^{2}\left(  \Omega\right)  \right)  ^{2}}^{2}\right)
^{1+M_{\ell}}\\
& \leq e^{\frac{s\mu_{1}\left(  1+M_{\ell}\right)  }{\left(  \ell+1\right)
h}}\left(  \left\Vert \left(  f_{1},f_{2}\right)  \left(  \cdot,T-\ell
h\right)  \right\Vert _{\left(  L^{2}\left(  \Omega\right)  \right)  ^{2}}%
^{2}\right)  ^{1+M_{\ell}}\\
& \leq K_{\ell}\,e^{\frac{s\mu_{1}\left(  1+M_{\ell}\right)  }{\left(
\ell+1\right)  h}}\left(  2\left\Vert \left(  f_{1},f_{2}\right)  \left(
\cdot,T\right)  \right\Vert _{\left(  L^{2}\left(  \Omega\right)  \right)
^{2}}^{2}\right)  \left(  2\left\Vert \left(  f_{1},f_{2}\right)  \left(
\cdot,T-2\ell h\right)  \right\Vert _{\left(  L^{2}\left(  \Omega\right)
\right)  ^{2}}^{2}\right)  ^{M_{\ell}}\\
& \leq2K_{\ell}\,e^{\frac{s\mu_{1}\left(  1+M_{\ell}\right)  }{\left(
\ell+1\right)  h}}\left(  2\left\Vert \left(  u_{1},u_{2}\right)  \left(
\cdot,0\right)  \right\Vert _{\left(  L^{2}\left(  \Omega\right)  \right)
^{2}}^{2}\right)  ^{M_{\ell}}\\
& \quad\times\left(  \left\Vert \left(  u_{1},u_{2}\right)  \left(
\cdot,T\right)  \right\Vert _{\left(  L^{2}\left(  B\left(  x_{0},r\right)
\right)  \right)  ^{2}}^{2}+e^{-\frac{s\mu_{0}}{h}}\left\Vert \left(
u_{1},u_{2}\right)  \left(  \cdot,0\right)  \right\Vert _{\left(  L^{2}\left(
\Omega\right)  \right)  ^{2}}^{2}\right)  \text{ .}%
\end{array}
\label{f5}%
\end{equation}
We now will choose $\ell>1$ sufficiently large to fulfill the inequality
$\frac{\mu_{1}\left(  1+M_{\ell}\right)  }{\left(  \ell+1\right)  }\leq
\frac{\mu_{0}}{2}$, so that $\frac{s\mu_{1}\left(  1+M_{\ell}\right)
}{\left(  \ell+1\right)  h}-\frac{s\mu_{0}}{h}\leq-\frac{s\mu_{0}}{2h}$. This
is possible since $\underset{\ell\rightarrow+\infty}{\text{lim}}\frac{\mu
_{1}\left(  1+M_{\ell}\right)  }{\left(  \ell+1\right)  }=0$ because (see
(\ref{defml})) $M_{\ell}\leq3e^{C_{1}}\frac{\left(  \ell+1\right)  ^{C_{0}}%
}{1-\left(  \frac{2}{3}\right)  ^{C_{0}}}$, with $C_{0}\in\left(  0,1\right)
$. Note that chosen in this way, $\ell$ depends on $|x_{0}|$, $r$ and
$d_{1},d_{2}$. \medskip

Thus,
\begin{equation}%
\begin{array}
[c]{ll}
& \quad\left(  \left\Vert \left(  u_{1},u_{2}\right)  \left(  \cdot,T\right)
\right\Vert _{\left(  L^{2}\left(  \Omega\right)  \right)  ^{2}}^{2}\right)
^{1+M_{\ell}}\\
& \leq2K_{\ell}\,\left(  2\left\Vert \left(  u_{1},u_{2}\right)  \left(
\cdot,0\right)  \right\Vert _{\left(  L^{2}\left(  \Omega\right)  \right)
^{2}}^{2}\right)  ^{M_{\ell}}\\
& \quad\times\left(  e^{\frac{s\mu_{0}}{2h}}\left\Vert \left(  u_{1}%
,u_{2}\right)  \left(  \cdot,T\right)  \right\Vert _{\left(  L^{2}\left(
B\left(  x_{0},r\right)  \right)  \right)  ^{2}}^{2}+e^{-\frac{s\mu_{0}}{2h}%
}\left\Vert \left(  u_{1},u_{2}\right)  \left(  \cdot,0\right)  \right\Vert
_{\left(  L^{2}\left(  \Omega\right)  \right)  ^{2}}^{2}\right)  \text{ .}%
\end{array}
\label{f6}%
\end{equation}
Since now $\ell>1$ is large but fixed, we denote $M:=M_{\ell}$ and
$K:=K_{\ell}$. We also take $s:=s_{2}$ only depending on $K_{0}$, $|x_{0}|$ and
$d_{1}$, $d_{2}$. Then, for any $h\in\left(  0,1\right]  $\ satisfying
$h<$min$\left(  1/\left(  2\ell\right)  ,T/\left(  4\ell\right)  \right)  $,
estimate (\ref{f6}) becomes
\begin{equation}%
\begin{array}
[c]{ll}
& \quad\left(  \left\Vert \left(  u_{1},u_{2}\right)  \left(  \cdot,T\right)
\right\Vert _{\left(  L^{2}\left(  \Omega\right)  \right)  ^{2}}^{2}\right)
^{1+M}\\
& \leq2^{1+M}K\left(  \left\Vert \left(  u_{1},u_{2}\right)  \left(
\cdot,0\right)  \right\Vert _{\left(  L^{2}\left(  \Omega\right)  \right)
^{2}}^{2}\right)  ^{M}\\
& \quad\times\left(  e^{\frac{s\mu_{0}}{2h}}\left\Vert \left(  u_{1}%
,u_{2}\right)  \left(  \cdot,T\right)  \right\Vert _{\left(  L^{2}\left(
B\left(  x_{0},r\right)  \right)  \right)  ^{2}}^{2}+e^{-\frac{s\mu_{0}}{2h}%
}\left\Vert \left(  u_{1},u_{2}\right)  \left(  \cdot,0\right)  \right\Vert
_{\left(  L^{2}\left(  \Omega\right)  \right)  ^{2}}^{2}\right)  \text{ .}%
\end{array}
\label{f7}%
\end{equation}
Recall that $1/\left(  2\ell\right)  \leq1$, therefore interpolation
inequality (\ref{f7}) holds for any $h>0$\ satisfying $h<$min$\left(
1/\left(  2\ell\right)  ,T/\left(  4\ell\right)  \right)  $.

On the other hand, for any $h\geq$min$\left(  1/\left(  2\ell\right)
,T/\left(  4\ell\right)  \right)  $, $1=e^{-\frac{s\mu_{0}}{2h}}e^{\frac
{s\mu_{0}}{2h}}\leq e^{-\frac{s\mu_{0}}{2h}}e^{\frac{s\mu_{0}}{2}\left(
2\ell+\frac{4\ell}{T}\right)  }$, which implies thanks to (\ref{3.3}):
\[
\left(  \left\Vert \left(  u_{1},u_{2}\right)  \left(  \cdot,T\right)
\right\Vert _{\left(  L^{2}\left(  \Omega\right)  \right)  ^{2}}^{2}\right)
^{1+M}\leq\left(  \left\Vert \left(  u_{1},u_{2}\right)  \left(
\cdot,0\right)  \right\Vert _{\left(  L^{2}\left(  \Omega\right)  \right)
^{2}}^{2}\right)  ^{1+M}e^{-\frac{s\mu_{0}}{2h}}e^{\frac{s\mu_{0}}{2}\left(
2\ell+\frac{4\ell}{T}\right)  }\text{ .}%
\]
Therefore, there exists $\mu_{2}$ depending only on $K_{0}$, $\left\vert
x_{0}\right\vert $, $r$ and $d_{1},d_{2}$, such that for any $h>0$\ satisfying
$h\geq$min$\left(  1/\left(  2\ell\right)  ,T/\left(  4\ell\right)  \right)
$, we have
\begin{equation}%
\begin{array}
[c]{ll}
& \quad\left(  \left\Vert \left(  u_{1},u_{2}\right)  \left(  \cdot,T\right)
\right\Vert _{\left(  L^{2}\left(  \Omega\right)  \right)  ^{2}}^{2}\right)
^{1+M}\\
& \leq e^{\mu_{2}\left(  1+\frac{1}{T}\right)  }\left(  \left\Vert \left(
u_{1},u_{2}\right)  \left(  \cdot,0\right)  \right\Vert _{\left(  L^{2}\left(
\Omega\right)  \right)  ^{2}}^{2}\right)  ^{M}\\
& \quad\times\left(  e^{\frac{s\mu_{0}}{2h}}\left\Vert \left(  u_{1}%
,u_{2}\right)  \left(  \cdot,T\right)  \right\Vert _{\left(  L^{2}\left(
B\left(  x_{0},r\right)  \right)  \right)  ^{2}}^{2}+e^{-\frac{s\mu_{0}}{2h}%
}\left\Vert \left(  u_{1},u_{2}\right)  \left(  \cdot,0\right)  \right\Vert
_{\left(  L^{2}\left(  \Omega\right)  \right)  ^{2}}^{2}\right)  \text{ .}%
\end{array}
\label{f8}%
\end{equation}

Consequently, one can conclude from (\ref{f7}) and (\ref{f8}) that there exists
$\mu_{3}$ depending only on $K_{0}$, $|x_{0}|$, $r$ and $d_{1},d_{2}$, such
that for any $h>0$,
\begin{equation}%
\begin{array}
[c]{ll}
& \quad\left(  \left\Vert \left(  u_{1},u_{2}\right)  \left(  \cdot,T\right)
\right\Vert _{\left(  L^{2}\left(  \Omega\right)  \right)  ^{2}}^{2}\right)
^{1+M}\\
& \leq e^{\mu_{3}\left(  1+\frac{1}{T}\right)  }\left(  \left\Vert \left(
u_{1},u_{2}\right)  \left(  \cdot,0\right)  \right\Vert _{\left(  L^{2}\left(
\Omega\right)  \right)  ^{2}}^{2}\right)  ^{M}\\
& \quad\times\left(  e^{\frac{s\mu_{0}}{2h}}\left\Vert \left(  u_{1}%
,u_{2}\right)  \left(  \cdot,T\right)  \right\Vert _{\left(  L^{2}\left(
B\left(  x_{0},r\right)  \right)  \right)  ^{2}}^{2}+e^{-\frac{s\mu_{0}}{2h}%
}\left\Vert \left(  u_{1},u_{2}\right)  \left(  \cdot,0\right)  \right\Vert
_{\left(  L^{2}\left(  \Omega\right)  \right)  ^{2}}^{2}\right)  \text{ .}%
\end{array}
\label{f9}%
\end{equation}
Now, one can minimize with respect to $h$ in $\left(  0,+\infty\right)  $ or
simply choose $h>0$ such that
\[
e^{-\frac{s\mu_{0}}{2h}}=\frac{\left\Vert \left(  u_{1},u_{2}\right)  \left(
\cdot,T\right)  \right\Vert _{\left(  L^{2}\left(  B\left(  x_{0},r\right)
\right)  \right)  ^{2}}}{\left\Vert \left(  u_{1},u_{2}\right)  \left(
\cdot,0\right)  \right\Vert _{\left(  L^{2}\left(  \Omega\right)  \right)
^{2}}}\text{ ,}%
\]
in order to obtain the desired estimate under the form
\[%
\begin{array}
[c]{ll}%
\left(  \left\Vert \left(  u_{1},u_{2}\right)  \left(  \cdot,T\right)
\right\Vert _{\left(  L^{2}\left(  \Omega\right)  \right)  ^{2}}^{2}\right)
^{1+\left(  1+2M\right)  } & \leq4e^{2\mu_{3}\left(  1+\frac{1}{T}\right)
}\left\Vert \left(  u_{1},u_{2}\right)  \left(  \cdot,T\right)  \right\Vert
_{\left(  L^{2}\left(  B\left(  x_{0},r\right)  \right)  \right)  ^{2}}^{2}\\
& \quad\times\left(  \left\Vert \left(  u_{1},u_{2}\right)  \left(
\cdot,0\right)  \right\Vert _{\left(  L^{2}\left(  \Omega\right)  \right)
^{2}}^{2}\right)  ^{\left(  1+2M\right)  }\text{ .}%
\end{array}
\]
This concludes the proof of Theorem \ref{mt2}.
\medskip

\section{Weight function and proof of Proposition \ref{prop2}}

\label{sec4}

We recall that
\[
\left\{
\begin{array}
[c]{ll}%
\varphi_{1}\left(  x\right)  =\psi\left(  x\right)  -\psi\left(  x_{0}\right)
\text{ ,} & \quad\text{for any}~x\in\overline{\Omega}\text{ ,}\\
\varphi_{3}\left(  x\right)  =-\psi\left(  x\right)  -\psi\left(
x_{0}\right)  \text{ ,} & \quad\text{for any}~x\in\overline{\Omega}\text{ ,}%
\end{array}
\right.
\]
and start this section with the

\begin{lem}
\label{lemtheta} Let $\vartheta=\left\{  \rho\leq\left\vert x\right\vert
<R\right\}  \subset\Omega$ be a neighborhood of $\partial\Omega$ with $\rho>0$
such that $x_{0}\notin\overline{\vartheta}$. Then there exists $c_{1}%
,c_{2},c_{3}>0$ depending only on $|x_{0}|$ and $\rho$, such that:

$\left(  i\right)  $ For any $x\in\Omega$,%
\[
\left\vert \nabla\varphi_{1}\left(  x\right)  \right\vert ^{2}\leq
c_{1}\left\vert \varphi_{1}\left(  x\right)  \right\vert \text{ and
}\left\vert \nabla\varphi_{3}\left(  x\right)  \right\vert ^{2}\leq
c_{1}\left\vert \varphi_{3}\left(  x\right)  \right\vert \text{ ;}%
\]

$\left(  ii\right)  $ For any $x\in\vartheta$,
\[
\left\vert \varphi_{1}\left(  x\right)  \right\vert \leq c_{2}\left\vert
\nabla\varphi_{1}\left(  x\right)  \right\vert ^{2}\text{ and }\left\vert
\varphi_{3}\left(  x\right)  \right\vert \leq c_{2}\left\vert \nabla
\varphi_{3}\left(  x\right)  \right\vert ^{2}\text{ ;}%
\]

$\left(  iii\right)  $ For any $x\in\Omega\left\backslash \vartheta\right.
$,
\[
\left\vert \varphi_{1}\left(  x\right)  \right\vert \leq c_{2}\left\vert
\nabla\varphi_{1}\left(  x\right)  \right\vert ^{2}\text{ and }\varphi
_{3}\left(  x\right)  -\varphi_{1}\left(  x\right)  \leq-c_{3}\text{ .}%
\]

\end{lem}

\bigskip

Proof of Lemma \ref{lemtheta}: Thanks to estimate (\ref{gpp}), and noticing
that $\left\vert \nabla\psi\left(  x\right)  \right\vert ^{2}=\left\vert
\nabla\varphi_{1}\left(  x\right)  \right\vert ^{2}$ and $\psi\left(
x_{0}\right)  -\psi\left(  x\right)  =\left\vert \varphi_{1}\left(  x\right)
\right\vert $, we conclude that $\left(  i\right)  -\left(  ii\right)  -\left(
iii\right)  $ holds for $\varphi_{1}$. For $\varphi_{3}$, we get $\left(
i\right)  -\left(  ii\right)  $ because $\left\vert \varphi_{3}\left(
x\right)  \right\vert >0$ for any $x\in\overline{\Omega}$ and $\left\vert
\nabla\varphi_{3}\left(  x\right)  \right\vert >0$ for any $x\in
\overline{\vartheta}$. Finally, $\varphi_{3}\left(  x\right)  -\varphi
_{1}\left(  x\right)  =-2\psi\left(  x\right)  <0$ for any $x\in
\Omega\left\backslash \vartheta\right.  $, which enables to complete the proof
of Lemma \ref{lemtheta}.

\medskip

We now turn to the \medskip

\textbf{{End of the proof of Proposition \ref{prop2}}}: First, we observe
that
\[
\eta_{i}=\frac{1}{2}\partial_{t}\Phi_{i}+\frac{1}{4}d_{i}\left\vert \nabla
\Phi_{i}\right\vert ^{2}=\frac{s}{\Gamma^{2}}\left(  -\frac{1}{2}\left\vert
\varphi_{i}\right\vert +\frac{1}{4}d_{i}s\left\vert \nabla\varphi
_{i}\right\vert ^{2}\right)  \text{ ,}%
\]
so that using Lemma \ref{lemtheta} $\left(  i\right)  $, there exists
$s_{0}\in(0,1]$ depending only on $d_{1}$, $d_{2}$, and $|x_{0}|$ such that
when $s\in(0,s_{0}]$,
\begin{equation}
\eta_{i}\left(  x,t\right)  \leq0\quad\text{ for any }\left(  x,t\right)
\in\Omega\times\left[  0,T\right]  \text{ .} \label{p21}%
\end{equation}
As a consequence, using the second part of identity (\ref{asds}), one can
deduce that for any $s\in(0,s_{0}]$, the estimate $\left\langle \mathcal{S}%
f,f\right\rangle \geq0$ holds.

\bigskip

Now, one can prove part $ii)$ of Proposition \ref{prop2}. Indeed, using
(\ref{3.1}) and the identities and estimates (\ref{prodo0}), $\Phi_{3}\leq
\Phi_{1}$, (\ref{pv3}),

we see that%
\[%
\begin{array}
[c]{ll}%
\left\Vert \digamma\right\Vert ^{2} & =\displaystyle\sum_{i=1,..,4}%
\displaystyle\int_{\Omega}\left\vert v_{i}\right\vert ^{2}e^{\Phi_{i}}\\
& \leq2\displaystyle\int_{\Omega}(|v_{1}|^{2}+|v_{2}|^{2})\,e^{\Phi_{1}}\\
& \leq2K_{0}\displaystyle\int_{\Omega}\bigg[|u_{1}|^{2}+|u_{2}|^{2}%
+(|u_{1}|^{2}+|u_{2}|^{2})^{2}\bigg]\,e^{\Phi_{1}}\\
& \leq4K_{0}\left(  \left\Vert f\right\Vert ^{2}+\displaystyle\int_{\Omega
}\left\vert u_{1}\right\vert ^{4}e^{\Phi_{1}}+\displaystyle\int_{\Omega
}\left\vert u_{2}\right\vert ^{4}e^{\Phi_{2}}\right)  \text{ .}%
\end{array}
\]
But thanks to H\"{o}lder and Sobolev inequalities (and denoting by $C_{Sob}$
the constant in this last inequality), for any $i\in\left\{  1,2\right\}$,
\[%
\begin{array}
[c]{ll}%
\displaystyle\int_{\Omega}\left\vert u_{i}\right\vert ^{4}e^{\Phi_{i}} &
\leq\left\Vert u_{i}^{2}\right\Vert _{L^{p}\left(  \Omega\right)  }\left\Vert
u_{i}^{2}e^{\Phi_{i}}\right\Vert _{L^{q}\left(  \Omega\right)  }\text{
whenever }\frac{1}{p}+\frac{1}{q}=1\\
& \leq\left(  \displaystyle\int_{\Omega}\left\vert u_{i}\right\vert
^{2p}\right)  ^{\frac{1}{p}}\left(  \displaystyle\int_{\Omega}\left\vert
f_{i}\right\vert ^{2q}\right)  ^{\frac{1}{q}}\text{ with }q=3\\
& \leq K_{0}\,C_{Sob}\left(  \displaystyle\int_{\Omega}\left\vert f_{i}%
\right\vert ^{2}+\displaystyle\int_{\Omega}\left\vert \nabla f_{i}\right\vert
^{2}\right)  \text{ ,}%
\end{array}
\]
using assumption (\ref{3.2}) in the last inequality, and remembering that
$H^{1}(\Omega)\subset L^{6}(\Omega)$ for n$\leq3$. Since $\eta_{i}\leq0$ and
therefore $d_{1}\int_{\Omega}\left\vert \nabla f_{1}\right\vert ^{2}+d_{2}%
\int_{\Omega}\left\vert \nabla f_{2}\right\vert ^{2}\leq\left\langle
\mathcal{S}f,f\right\rangle $, this gives the desired estimate when $C_{1}
\geq 4K_{0}\,$max$\left(  1+K_{0}C_{Sob},(\underset{i}{\text{min}}d_{i}
)^{-1}K_{0}C_{Sob}\right)  $.
\medskip

It remains to prove part $iii)$ of Proposition \ref{prop2}.
\medskip

We recall that
\[
\left\langle \mathcal{S}^{\prime}f,f\right\rangle =\sum_{i=1,..,4}\int
_{\Omega}\left(  -\partial_{t}\eta_{i}\right)  \left\vert f_{i}\right\vert
^{2}\text{ .}%
\]
Moreover, using the definition of $\mathcal{S}f$ and $\mathcal{A}f$, the
bracket $2\left\langle \mathcal{S}f,\mathcal{A}f\right\rangle $ writes
\[
2\left\langle \mathcal{S}f,\mathcal{A}f\right\rangle =2\sum_{i=1,..,4}%
\int_{\Omega}\left(  d_{i}\Delta f_{i}+\eta_{i}f_{i}\right)  \left(
d_{i}\nabla\Phi_{i}\cdot\nabla f_{i}+\frac{1}{2}d_{i}\Delta\Phi_{i}%
f_{i}\right)  \text{ .}%
\]
Then, four integrations by parts give%
\[%
\begin{array}
[c]{ll}%
2\left\langle \mathcal{S}f,\mathcal{A}f\right\rangle  & =\displaystyle\sum
_{i=1,..,4}\left(  -2d_{i}^{2}\displaystyle\int_{\Omega}\nabla f_{i}\nabla
^{2}\Phi_{i}\nabla f_{i}-d_{i}^{2}\displaystyle\int_{\Omega}\nabla f_{i}%
\Delta\nabla\Phi_{i}f_{i}-d_{i}\displaystyle\int_{\Omega}\nabla\eta_{i}%
\cdot\nabla\Phi_{i}\left\vert f_{i}\right\vert ^{2}\right) \\
& \quad+\displaystyle2\sum_{i=1,..,4}d_{i}^{2}\displaystyle\int_{\partial
\Omega}\partial_{n}f_{i}\nabla\Phi_{i}\cdot\nabla f_{i}-\displaystyle\sum
_{i=1,..,4}d_{i}^{2}\int_{\partial\Omega}\partial_{n}\Phi_{i}\left\vert \nabla
f_{i}\right\vert ^{2}\\
& \quad+\displaystyle\sum_{i=1,..,4}d_{i}^{2}\displaystyle\int_{\partial
\Omega}\partial_{n}f_{i}\Delta\Phi_{i}f_{i}+\displaystyle\sum_{i=1,..,4}%
\displaystyle d_{i}\int_{\partial\Omega}\eta_{i}\partial_{n}\Phi_{i}\left\vert
f_{i}\right\vert ^{2}\text{ .}%
\end{array}
\]
Indeed,
\[
\int_{\Omega}\Delta f_{i}\nabla\Phi_{i}\cdot\nabla f_{i}=\int_{\partial\Omega
}\partial_{n}f_{i}\nabla\Phi_{i}\cdot\nabla f_{i}-\int_{\Omega}\nabla
f_{i}\nabla^{2}\Phi_{i}\nabla f_{i}-\int_{\Omega}\nabla f_{i}\nabla^{2}%
f_{i}\nabla\Phi_{i}\text{ ,}%
\]
and
\[
\int_{\Omega}\nabla f_{i}\nabla^{2}f_{i}\nabla\Phi_{i}=\frac{1}{2}%
\int_{\partial\Omega}\partial_{n}\Phi_{i}\left\vert \nabla f_{i}\right\vert
^{2}-\frac{1}{2}\int_{\Omega}\Delta\Phi_{i}\left\vert \nabla f_{i}\right\vert
^{2}\text{ .}%
\]
Second,
\[
\int_{\Omega}\Delta f_{i}\Delta\Phi_{i}f_{i}=\int_{\partial\Omega}\partial
_{n}f_{i}\Delta\Phi_{i}f_{i}-\int_{\Omega}\nabla f_{i}\Delta\nabla\Phi
_{i}f_{i}-\int_{\Omega}\Delta\Phi_{i}\left\vert \nabla f_{i}\right\vert
^{2}\text{ .}%
\]
Third,
\[
2\int_{\Omega}\eta_{i}f_{i}\nabla\Phi_{i}\cdot\nabla f_{i}=\int_{\partial
\Omega}\eta_{i}\partial_{n}\Phi_{i}\left\vert f_{i}\right\vert ^{2}%
-\int_{\Omega}\nabla\eta_{i}\cdot\nabla\Phi_{i}\left\vert f_{i}\right\vert
^{2}-\int_{\Omega}\eta_{i}\Delta\Phi_{i}\left\vert f_{i}\right\vert ^{2}\text{
.}%
\]

Then, we see that
\[%
\begin{array}
[c]{ll}%
\left\langle \mathcal{S}^{\prime}f,f\right\rangle +2\left\langle
\mathcal{S}f,\mathcal{A}f\right\rangle  & =\displaystyle\sum_{i=1,..,4}\left(
-2d_{i}^{2}\displaystyle\int_{\Omega}\nabla f_{i}\nabla^{2}\Phi_{i}\nabla
f_{i}-d_{i}^{2}\int_{\Omega}\nabla f_{i}\Delta\nabla\Phi_{i}f_{i}\right) \\
& \quad+\displaystyle\sum_{i=1,..,4}\displaystyle\int_{\Omega}\left(
-\partial_{t}\eta_{i}-d_{i}\nabla\eta_{i}\cdot\nabla\Phi_{i}\right)
\left\vert f_{i}\right\vert ^{2}\\
& \quad+\displaystyle2\sum_{i=1,..,4}d_{i}^{2}\displaystyle\int_{\partial
\Omega}\partial_{n}f_{i}\nabla\Phi_{i}\cdot\nabla f_{i}-\displaystyle\sum
_{i=1,..,4}d_{i}^{2}\int_{\partial\Omega}\partial_{n}\Phi_{i}\left\vert \nabla
f_{i}\right\vert ^{2}\\
& \quad+\displaystyle\sum_{i=1,..,4}d_{i}^{2}\displaystyle\int_{\partial
\Omega}\partial_{n}f_{i}\Delta\Phi_{i}f_{i}+\displaystyle\sum_{i=1,..,4}%
\displaystyle d_{i}\int_{\partial\Omega}\eta_{i}\partial_{n}\Phi_{i}\left\vert
f_{i}\right\vert ^{2}\text{ .}%
\end{array}
\]
Next, the computation of $\partial_{t}\eta_{i}+d_{i}\nabla\eta_{i}\cdot
\nabla\Phi_{i}$ gives%
\[%
\begin{array}
[c]{ll}%
\partial_{t}\eta_{i}+d_{i}\nabla\eta_{i}\cdot\nabla\Phi_{i} &
=\displaystyle\frac{1}{2}\partial_{t}^{2}\Phi_{i}+d_{i}\nabla\partial_{t}%
\Phi_{i}\cdot\nabla\Phi_{i}+\displaystyle\frac{1}{2}d_{i}^{2}\nabla\Phi
_{i}\nabla^{2}\Phi_{i}\nabla\Phi_{i}\\
& =\displaystyle\frac{1}{\Gamma}\partial_{t}\Phi_{i}+\displaystyle\frac
{1}{\Gamma}d_{i}\left\vert \nabla\Phi_{i}\right\vert ^{2}+\displaystyle\frac
{s}{2\Gamma}d_{i}^{2}\nabla\Phi_{i}\nabla^{2}\varphi_{i}\nabla\Phi_{i}\\
& =\displaystyle\frac{2}{\Gamma}\left(  \frac{1}{2}\partial_{t}\Phi_{i}%
+\frac{1}{4}d_{i}\left\vert \nabla\Phi_{i}\right\vert ^{2}\right)
+\displaystyle\frac{1}{2\Gamma}d_{i}\left\vert \nabla\Phi_{i}\right\vert
^{2}+\displaystyle\frac{s}{2\Gamma}d_{i}^{2}\nabla\Phi_{i}\nabla^{2}%
\varphi_{i}\nabla\Phi_{i},
\end{array}
\]
since $\partial_{t}^{2}\Phi_{i}=\frac{2}{\Gamma}\partial_{t}\Phi_{i}$ and
$\partial_{t}\Phi_{i}=\frac{1}{\Gamma}\Phi_{i}$. Therefore, by definition
(\ref{eta}) of $\eta_{i}$,
\[
-\partial_{t}\eta_{i}-d_{i}\nabla\eta_{i}\cdot\nabla\Phi_{i}=-\frac{2}{\Gamma
}\eta_{i}-\frac{1}{2\Gamma}d_{i}\left\vert \nabla\Phi_{i}\right\vert
^{2}-\frac{s}{2\Gamma}d_{i}^{2}\nabla\Phi_{i}\nabla^{2}\varphi_{i}\nabla
\Phi_{i}\,\text{,}%
\]
and one can conclude that
\begin{equation}%
\begin{array}
[c]{ll}%
\left\langle \mathcal{S}^{\prime}f,f\right\rangle +2\left\langle
\mathcal{S}f,\mathcal{A}f\right\rangle  & =\displaystyle\sum_{i=1,..,4}\left(
-2d_{i}^{2}\displaystyle\int_{\Omega}\nabla f_{i}\nabla^{2}\Phi_{i}\nabla
f_{i}-d_{i}^{2}\int_{\Omega}\nabla f_{i}\Delta\nabla\Phi_{i}f_{i}\right) \\
& \quad+\displaystyle\frac{1}{\Gamma}\sum_{i=1,..,4}\displaystyle\int_{\Omega
}\left(  -2\eta_{i}-\frac{1}{2}d_{i}\left\vert \nabla\Phi_{i}\right\vert
^{2}-\frac{s}{2}d_{i}^{2}\nabla\Phi_{i}\nabla^{2}\varphi_{i}\nabla\Phi
_{i}\right)  \left\vert f_{i}\right\vert ^{2}\\
& \quad+\displaystyle2\sum_{i=1,..,4}d_{i}^{2}\displaystyle\int_{\partial
\Omega}\partial_{n}f_{i}\nabla\Phi_{i}\cdot\nabla f_{i}-\displaystyle\sum
_{i=1,..,4}d_{i}^{2}\int_{\partial\Omega}\partial_{n}\Phi_{i}\left\vert \nabla
f_{i}\right\vert ^{2}\\
& \quad+\displaystyle\sum_{i=1,..,4}d_{i}^{2}\displaystyle\int_{\partial
\Omega}\partial_{n}f_{i}\Delta\Phi_{i}f_{i}+\displaystyle\sum_{i=1,..,4}%
\displaystyle d_{i}\int_{\partial\Omega}\eta_{i}\partial_{n}\Phi_{i}\left\vert
f_{i}\right\vert ^{2}\text{ .}%
\end{array}
\label{4.1}%
\end{equation}

First we estimate the contribution of the gradient terms. We observe that
(since $\psi$ is smooth on $\overline{\Omega}$),
\[
\left\vert d_{i}\nabla^{2}\Phi_{i}\right\vert +\left\vert d_{i}\Delta
\nabla\Phi_{i}\right\vert \leq\frac{C_{2}\,s}{\Gamma}\text{ ,}%
\]
where $C_{2}$ only depends on $d_{1}$, $d_{2}$, and $|x_{0}|$.

As a consequence, using Young's inequality,
\begin{equation}%
\begin{array}
[c]{ll}%
\displaystyle\sum_{i=1,..,4}\left(  -2d_{i}^{2}\displaystyle\int_{\Omega
}\nabla f_{i}\nabla^{2}\Phi_{i}\nabla f_{i}-d_{i}^{2}\displaystyle\int
_{\Omega}\nabla f_{i}\Delta\nabla\Phi_{i}f_{i}\right)  & \leq
\displaystyle\frac{C_{3}s}{\Gamma}\sum_{i=1,..,4}d_{i}\int_{\Omega}\left\vert
\nabla f_{i}\right\vert ^{2}+\frac{C_{3}s}{\Gamma}\left\Vert f\right\Vert
^{2}\\
& \leq\displaystyle\frac{C_{3}s}{\Gamma}\sum_{i=1,..,4}d_{i}\int_{\Omega
}\left\vert \nabla f_{i}\right\vert ^{2}+\displaystyle\frac{C_{3}}%
{h}\left\Vert f\right\Vert ^{2}\text{ ,}%
\end{array}
\label{4.2}%
\end{equation}
where we used in the last line the inequalities $\frac{1}{\Gamma}\leq\frac
{1}{h}$ and $s\in\left(  0,1\right]  $. Here, $C_{3}$ only depends on $d_{1}$,
$d_{2}$, and $|x_{0}|$.

\begin{lem}
\label{lemnn} There exists a constant $C_{5}>0$ depending only on $|x_{0}|$
and $d_{1}$, $d_{2}$ such that for any $s\in\left(  0,1\right]  $,
\[%
\begin{array}
[c]{ll}
& \quad\displaystyle2\sum_{i=1,..,4}d_{i}^{2}\displaystyle\int_{\partial
\Omega}\partial_{n}f_{i}\nabla\Phi_{i}\cdot\nabla f_{i}-\displaystyle\sum
_{i=1,..,4}d_{i}^{2}\displaystyle\int_{\partial\Omega}\partial_{n}\Phi
_{i}\left\vert \nabla f_{i}\right\vert ^{2}\\
& \quad+\displaystyle\sum_{i=1,..,4}d_{i}^{2}\displaystyle\int_{\partial
\Omega}\partial_{n}f_{i}\Delta\Phi_{i}f_{i}+\displaystyle\sum_{i=1,..,4}%
\displaystyle d_{i}\int_{\partial\Omega}\eta_{i}\partial_{n}\Phi_{i}\left\vert
f_{i}\right\vert ^{2}\\
& \leq\displaystyle\frac{C_{5}s}{\Gamma}\displaystyle\sum_{i=1,..,4}%
d_{i}\displaystyle\int_{\Omega}\left\vert \nabla f_{i}\right\vert
^{2}+\displaystyle\frac{C_{5}s}{\Gamma}\displaystyle\sum_{i=1,2}%
d_{i}\displaystyle\int_{\Omega}\left\vert \nabla\Phi_{i}\right\vert
^{2}\left\vert f_{i}\right\vert ^{2}+\displaystyle\frac{C_{5}}{h^{2}%
}\left\Vert f\right\Vert ^{2}\text{ .}%
\end{array}
\]

\end{lem}

\noindent\textbf{{Proof of Lemma \ref{lemnn}}}: We claim that
$\displaystyle\sum_{i=1,..,4}d_{i}\displaystyle\int_{\partial\Omega}\eta
_{i}\partial_{n}\Phi_{i}\left\vert f_{i}\right\vert ^{2}=0$.

We first observe that thanks to identities (\ref{prodo3}), (\ref{prodo4}), the
following extra identities hold:
\begin{equation}
\partial_{t}\Phi_{1}=\partial_{t}\Phi_{3}\text{ ,}\qquad\left\vert \nabla
\Phi_{1}\right\vert =\left\vert \nabla\Phi_{3}\right\vert =0\quad
{\hbox{ on }}\quad\partial\Omega\times\left(  0,T\right)  \text{
.}\label{prodon}%
\end{equation}
Then, since $\eta_{i}=\frac{1}{2}\partial_{t}\Phi_{i}+\frac{1}{4}%
d_{i}\left\vert \nabla\Phi_{i}\right\vert ^{2}$,
\[%
\begin{array}
[c]{ll}%
\displaystyle\sum_{i=1,..,4}\displaystyle\int_{\partial\Omega}d_{i}\eta
_{i}\partial_{n}\Phi_{i}\left\vert f_{i}\right\vert ^{2} & =\displaystyle\sum
_{i=1,..,4}\displaystyle\int_{\partial\Omega}d_{i}\left(  \frac{1}{2}%
\partial_{t}\Phi_{i}+\frac{1}{4}d_{i}\left\vert \nabla\Phi_{i}\right\vert
^{2}\right)  \partial_{n}\Phi_{i}\left\vert u_{i}\right\vert ^{2}e^{\Phi_{i}%
}\\
& =\displaystyle\int_{\partial\Omega}d_{1}\left(  \frac{1}{2}\partial_{t}%
\Phi_{1}+\frac{1}{4}d_{1}\left\vert \nabla\Phi_{1}\right\vert ^{2}\right)
\partial_{n}\Phi_{1}\left\vert u_{1}\right\vert ^{2}e^{\Phi_{1}}\\
& \quad+\displaystyle\int_{\partial\Omega}d_{2}\left(  \frac{1}{2}\partial
_{t}\Phi_{1}+\frac{1}{4}d_{2}\left\vert \nabla\Phi_{1}\right\vert ^{2}\right)
\partial_{n}\Phi_{1}\left\vert u_{2}\right\vert ^{2}e^{\Phi_{1}}\\
& \quad+\displaystyle\int_{\partial\Omega}d_{1}\left(  \frac{1}{2}\partial
_{t}\Phi_{3}+\frac{1}{4}d_{1}\left\vert \nabla\Phi_{3}\right\vert ^{2}\right)
\partial_{n}\Phi_{3}\left\vert u_{1}\right\vert ^{2}e^{\Phi_{3}}\\
& \quad+\displaystyle\int_{\partial\Omega}d_{2}\left(  \frac{1}{2}\partial
_{t}\Phi_{3}+\frac{1}{4}d_{2}\left\vert \nabla\Phi_{3}\right\vert ^{2}\right)
\partial_{n}\Phi_{3}\left\vert u_{2}\right\vert ^{2}e^{\Phi_{3}}\text{ ,}%
\end{array}
\]
where in the second line, we used identities (\ref{prodo0}), (\ref{prodo3}),
(\ref{prodo4}), (\ref{prodon}), (\ref{dd1}). This completes the claim.

\bigskip

We then observe that
\begin{equation}
\label{iid}2\sum_{i=1,..,4}d_{i}^{2}\int_{\partial\Omega}\partial_{n}%
f_{i}\nabla\Phi_{i}\cdot\nabla f_{i}-\sum_{i=1,..,4}d_{i}^{2}\int
_{\partial\Omega} \partial_{n}\Phi_{i}\left\vert \nabla f_{i}\right\vert
^{2}=0\text{ .}%
\end{equation}
Indeed, since $\nabla\Phi_{i}=\partial_{n}\Phi_{i}\overrightarrow{n}$ on
$\partial\Omega\times\left(  0,T\right)  $, we see first that
\[%
\begin{array}
[c]{ll}%
2\displaystyle\sum_{i=1,..,4}d_{i}^{2}\displaystyle\int_{\partial\Omega
}\partial_{n}f_{i}\nabla\Phi_{i}\cdot\nabla f_{i} & =2\displaystyle\sum
_{i=1,..,4}d_{i}^{2}\displaystyle\int_{\partial\Omega}\partial_{n}\Phi
_{i}\left\vert \partial_{n}f_{i}\right\vert ^{2}\\
& =2\displaystyle\sum_{i=1,..,4}d_{i}^{2}\displaystyle\int_{\partial\Omega
}\partial_{n}\Phi_{i}\left\vert \frac{1}{2}\partial_{n}\Phi_{i}f_{i}%
\right\vert ^{2}\text{ because }\partial_{n}f_{i}=\frac{1}{2}\partial_{n}%
\Phi_{i}f_{i}\\
& =0 \text{ ,}%
\end{array}
\]
thanks to identities (\ref{prodo0}), (\ref{prodo3}), (\ref{prodo4}).

We then observe that on $\partial\Omega\times\left(  0,T\right)  $,
\[%
\begin{array}
[c]{ll}%
\left\vert \nabla f_{i}\right\vert ^{2} & =\left\vert \nabla u_{i}e^{\Phi
_{i}/2}+u_{i}\frac{1}{2}\nabla\Phi_{i}e^{\Phi_{i}/2}\right\vert ^{2}\\
& =\left\vert \partial_{\tau}u_{i}\overrightarrow{\tau}+u_{i}\frac{1}%
{2}\partial_{n}\Phi_{i}\overrightarrow{n}\right\vert ^{2}e^{\Phi_{i}}\text{
because }\partial_{n}u_{i}=0\text{ and }\Phi_{i}\left\vert _{\partial\Omega
}\right.  =\text{constant}\\
& =\left(  \left\vert \partial_{\tau}u_{i}\right\vert ^{2}+\left\vert \frac
{1}{2}u_{i}\partial_{n}\Phi_{i}\right\vert ^{2}\right)  e^{\Phi_{i}}\text{ .}%
\end{array}
\]
As a consequence, $-\displaystyle\sum_{i=1,..,4}d_{i}^{2}\displaystyle\int
_{\partial\Omega}\partial_{n}\Phi_{i}\left\vert \nabla f_{i}\right\vert
^{2}=0$, where we used identities (\ref{prodo3}), (\ref{prodo4}).

We complete in this way the proof of identity (\ref{iid}). \medskip

Next, it remains to treat the contribution of $\displaystyle\sum
_{i=1,..,4}d_{i}^{2}\displaystyle\int_{\partial\Omega}\partial_{n}f_{i}%
\Delta\Phi_{i}f_{i}$. We introduce $C_{4}:=\underset{\overline{\Omega}%
}{\text{max}}\left\vert \Delta\psi\right\vert $. Note that $C_{4}$ only
depends on $|x_{0}|$.
\[%
\begin{array}
[c]{ll}%
\displaystyle\sum_{i=1,..,4}d_{i}^{2}\displaystyle\int_{\partial\Omega
}\partial_{n}f_{i}\Delta\Phi_{i}f_{i} & =\displaystyle\sum_{i=1,..,4}d_{i}%
^{2}\displaystyle\int_{\partial\Omega}\frac{1}{2}\partial_{n}\Phi_{i}%
\Delta\Phi_{i}\left\vert f_{i}\right\vert ^{2}\text{ because }\partial
_{n}f_{i}=\frac{1}{2}\partial_{n}\Phi_{i}f_{i}\\
& \leq\displaystyle\frac{C_{4}s}{2\Gamma}\sum_{i=1,..,4}d_{i}^{2}%
\displaystyle\int_{\partial\Omega}\left\vert \partial_{n}\Phi_{i}\right\vert
\left\vert f_{i}\right\vert ^{2}\text{ because }\left\vert \Delta\Phi
_{i}\right\vert =\frac{s}{\Gamma}\left\vert \Delta\varphi_{i}\right\vert
\leq\frac{C_{4}s}{\Gamma}\\
& \leq\displaystyle\frac{C_{4}s}{\Gamma}\sum_{i=1,2}d_{i}^{2}\displaystyle\int
_{\partial\Omega}\left\vert \partial_{n}\Phi_{i}\right\vert \left\vert
f_{i}\right\vert ^{2}\\
& =\displaystyle\frac{C_{4}s}{\Gamma}\sum_{i=1,2}d_{i}^{2}\displaystyle\int
_{\partial\Omega}\left(  -\partial_{n}\Phi_{i}\right)  \left\vert
f_{i}\right\vert ^{2}\text{ because }\partial_{n}\Phi_{1}\leq0\text{ and }%
\Phi_{1}=\Phi_{2}\text{ .}%
\end{array}
\]
In the third line we used $d_{3}^{2}\displaystyle\int_{\partial\Omega
}\left\vert \partial_{n}\Phi_{3}\right\vert \left\vert f_{3}\right\vert
^{2}=d_{1}^{2}\displaystyle\int_{\partial\Omega}\left\vert \partial_{n}%
\Phi_{1}\right\vert \left\vert f_{1}\right\vert ^{2}$, which holds thanks to
identities (\ref{dd1}), (\ref{prodo4}), and $f_{3}=f_{1}$ on $\partial
\Omega\times\left(  0,T\right)  $. Similar computations hold for $d_{4}%
^{2}\displaystyle\int_{\partial\Omega}\left\vert \partial_{n}\Phi
_{4}\right\vert \left\vert f_{4}\right\vert ^{2}$.

Then, thanks to an integration by parts,
\[%
\begin{array}
[c]{ll}%
\displaystyle\int_{\partial\Omega}\left(  -\partial_{n}\Phi_{i}\right)
\left\vert f_{i}\right\vert ^{2} & =-2\displaystyle\int_{\Omega}\nabla
f_{i}\cdot\nabla\Phi_{i}f_{i}-\displaystyle\int_{\Omega}\Delta\Phi
_{i}\left\vert f_{i}\right\vert ^{2}\\
& \leq\displaystyle\int_{\Omega}\left\vert \nabla f_{i}\right\vert
^{2}+\displaystyle\int_{\Omega}\left\vert \nabla\Phi_{i}\right\vert
^{2}\left\vert f_{i}\right\vert ^{2}+\displaystyle\frac{C_{4}s}{h}\left\Vert
f\right\Vert ^{2}\text{ ,}%
\end{array}
\]
using Young's inequality and the estimate $\left\vert \Delta\Phi
_{i}\right\vert =\frac{s}{\Gamma}\left\vert \Delta\varphi_{i}\right\vert
\leq\frac{C_{4}s}{h}$. Therefore, one can conclude that for any $s\in\left(
0,1\right]  $,%
\[%
\begin{array}
[c]{ll}%
\displaystyle\sum_{i=1,..,4}d_{i}^{2}\displaystyle\int_{\partial\Omega
}\partial_{n}f_{i}\Delta\Phi_{i}f_{i} & \leq\underset{i}{\text{max}}\,
d_{i}\displaystyle\frac{C_{4}\,s}{\Gamma}\displaystyle\sum_{i=1,..,4}%
d_{i}\displaystyle\int_{\Omega}\left\vert \nabla f_{i}\right\vert ^{2}\\
& \quad+\,\underset{i}{\text{max}}\,d_{i}\displaystyle\frac{C_{4}\,s}{\Gamma
}\displaystyle\sum_{i=1,2}d_{i}\displaystyle\int_{\Omega}\left\vert \nabla
\Phi_{i}\right\vert ^{2}\left\vert f_{i}\right\vert ^{2}+(\underset
{i}{\text{max}}\,d_{i})^{2}\displaystyle\frac{C_{4}^{2}}{h^{2}}\left\Vert
f\right\Vert ^{2}\text{ .}%
\end{array}
\]
This completes the proof of Lemma \ref{lemnn}. \medskip

Finally, we estimate the contribution of
\[
\frac{1}{\Gamma}\int_{\Omega}\left(  -2\eta_{i}-\frac{1}{2}d_{i}\left\vert
\nabla\Phi_{i}\right\vert ^{2}-\frac{s}{2}d_{i}^{2}\nabla\Phi_{i}\nabla
^{2}\varphi_{i}\nabla\Phi_{i}\right)  \left\vert f_{i}\right\vert ^{2}%
+\frac{C_{5}s}{\Gamma}d_{i}\int_{\Omega}\left\vert \nabla\Phi_{i}\right\vert
^{2}\left\vert f_{i}\right\vert ^{2}\text{ ,}%
\]
where $C_{5}$ is the constant appearing in Lemma \ref{lemnn}.

\bigskip

\begin{lem}
\label{prop4} There exists $s_{1}\in(0,1]$ and $C_{7}>0$, both depending only
on $|x_{0}|$ and $d_{1}$, $d_{2}$ such that when $s\in\left(  0,s_{1}\right]
$,
\[%
\begin{array}
[c]{ll}
& \quad\displaystyle\frac{1}{\Gamma}\sum_{i=1,..,4}\displaystyle\int_{\Omega
}\left(  -2\eta_{i}-\frac{1}{2}d_{i}\left\vert \nabla\Phi_{i}\right\vert
^{2}-\frac{s}{2}d_{i}^{2}\nabla\Phi_{i}\nabla^{2}\varphi_{i}\nabla\Phi
_{i}\right)  \left\vert f_{i}\right\vert ^{2}+\displaystyle\frac{C_{5}%
s}{\Gamma}\sum_{i=1,2}d_{i}\displaystyle\int_{\Omega}\left\vert \nabla\Phi
_{i}\right\vert ^{2}\left\vert f_{i}\right\vert ^{2}\\
& \leq\displaystyle\frac{1}{\Gamma}\sum_{i=1,..,4}\left(  2-\frac{1}{4}%
d_{i}\frac{s}{c_{2}}\right)  \displaystyle\int_{\Omega}\left(  -\eta
_{i}\right)  \left\vert f_{i}\right\vert ^{2}+\displaystyle\frac{C_{7}}%
{h}\left\Vert f\right\Vert ^{2}\text{ .}%
\end{array}
\]

\end{lem}

\noindent\textbf{{Proof of Lemma \ref{prop4}:}} First observe that
\[%
\begin{array}
[c]{ll}
& \quad\left(  -2\eta_{i}-\displaystyle\frac{1}{2}d_{i}\left\vert \nabla
\Phi_{i}\right\vert ^{2}-\frac{s}{2}d_{i}^{2}\nabla\Phi_{i}\nabla^{2}%
\varphi_{i}\nabla\Phi_{i}\right)  +C_{5}sd_{i}\left\vert \nabla\Phi
_{i}\right\vert ^{2}\\
& \leq-2\eta_{i}+\left(  -\displaystyle\frac{1}{2}+s\left(  \displaystyle\frac
{1}{2}\underset{i}{\text{max}}\,d_{i}\,\underset{\overline{\Omega}}{\text{max}%
}\left\vert \nabla^{2}\psi\right\vert +C_{5}\right)  \right)  d_{i}\left\vert
\nabla\Phi_{i}\right\vert ^{2}\\
& \leq-\left(  2\eta_{i}+\displaystyle\frac{1}{8}d_{i}\left\vert \nabla
\Phi_{i}\right\vert ^{2}\right)
\end{array}
\]
for any $s\in\left(  0,s_{1}\right]  $ if $s_{1}>0$ is well chosen. Indeed, it
is sufficient to take $s_{1}:=\frac{3}{8}(\frac{1}{2}\underset{i}{\text{max}%
}\,d_{i}\,\underset{\overline{\Omega}}{\text{max}}\left\vert \nabla^{2}%
\psi\right\vert +C_{5})^{-1}$.

Now, from
\[%
\begin{array}
[c]{ll}
& \quad\displaystyle\sum_{i=1,..,4}\displaystyle\int_{\Omega}\left(
-2\eta_{i}-\frac{1}{2}d_{i}\left\vert \nabla\Phi_{i}\right\vert ^{2}-\frac
{s}{2}d_{i}^{2}\nabla\Phi_{i}\nabla^{2}\varphi_{i}\nabla\Phi_{i}\right)
\left\vert f_{i}\right\vert ^{2}+C_{5}\,s\sum_{i=1,2}d_{i}\displaystyle\int
_{\Omega}\left\vert \nabla\Phi_{i}\right\vert ^{2}\left\vert f_{i}\right\vert
^{2}\\
& \leq\displaystyle\sum_{i=1,..,4}\left(  2\int_{\Omega}\left(  -\eta
_{i}\right)  \left\vert f_{i}\right\vert ^{2}-\frac{1}{8}d_{i}\int_{\Omega
}\left\vert \nabla\Phi_{i}\right\vert ^{2}\left\vert f_{i}\right\vert
^{2}\right)  \text{ ,}%
\end{array}
\]
we want to achieve
\[%
\begin{array}
[c]{ll}
& \quad\displaystyle\sum_{i=1,..,4}\displaystyle\int_{\Omega}\left(
-2\eta_{i}-\displaystyle\frac{1}{2}d_{i}\left\vert \nabla\Phi_{i}\right\vert
^{2}-\displaystyle\frac{s}{2}d_{i}^{2}\nabla\Phi_{i}\nabla^{2}\varphi
_{i}\nabla\Phi_{i}\right)  \left\vert f_{i}\right\vert ^{2}+C_{5}\text{
}s\displaystyle\sum_{i=1,2}d_{i}\displaystyle\int_{\Omega}\left\vert
\nabla\Phi_{i}\right\vert ^{2}\left\vert f_{i}\right\vert ^{2}\\
& \leq\displaystyle\sum_{i=1,..,4}\left(  2-\displaystyle\frac{1}{4}d_{i}%
\frac{s}{c_{2}}\right)  \displaystyle\int_{\Omega}\left(  -\eta_{i}\right)
\left\vert f_{i}\right\vert ^{2}+C_{7}\left\Vert f\right\Vert ^{2}\text{ .}%
\end{array}
\]
We will treat separately the case $i=1$ (which is similar to the case $i=2$
since $\Phi_{2}=\Phi_{1}$) and the case $i=3$ (which is similar to the case
$i=4$ since $\Phi_{4}=\Phi_{3}$).

\bigskip

Thanks to Lemma \ref{lemtheta} $\left(  ii\right)  -\left(  iii\right)  $, for
any $x\in\Omega$, $\left\vert \varphi_{1}\left(  x\right)  \right\vert \leq
c_{2}\left\vert \nabla\varphi_{1}\left(  x\right)  \right\vert ^{2}$. This
implies
\[
-\left\vert \nabla\Phi_{1}\right\vert ^{2}=-\frac{s^{2}}{\Gamma^{2}}\left\vert
\nabla\varphi_{1}\right\vert ^{2}\leq-\frac{s^{2}}{c_{2}\Gamma^{2}}\left\vert
\varphi_{1}\right\vert =\frac{2s}{c_{2}}\left(  -\frac{s}{2\Gamma^{2}%
}\left\vert \varphi_{1}\right\vert \right)  \leq\frac{2s}{c_{2}}\eta_{1}\text{
.}%
\]
Therefore,
\[
-\frac{1}{8}d_{1}\int_{\Omega}\left\vert \nabla\Phi_{1}\right\vert
^{2}\left\vert f_{1}\right\vert ^{2}\leq\frac{1}{4}d_{1}\frac{s}{c_{2}}%
\int_{\Omega}\eta_{1}\left\vert f_{1}\right\vert ^{2}%
\]
and similarly for $i=2$, one has $-\displaystyle\frac{1}{8}d_{2}\int_{\Omega
}\left\vert \nabla\Phi_{2}\right\vert ^{2}\left\vert f_{2}\right\vert ^{2}%
\leq\frac{1}{4}d_{2}\displaystyle\frac{s}{c_{2}}\int_{\Omega}\eta
_{2}\left\vert f_{2}\right\vert ^{2}$. Consequently,
\[
\sum_{i=1,2}\left(  -\frac{1}{8}d_{i}\int_{\Omega}\left\vert \nabla\Phi
_{i}\right\vert ^{2}\left\vert f_{i}\right\vert ^{2}\right)  \leq\sum
_{i=1,2}\frac{1}{4}d_{i}\frac{s}{c_{2}}\int_{\Omega}\eta_{i}\left\vert
f_{i}\right\vert ^{2}=\sum_{i=1,2}\left(  -\frac{1}{4}d_{i}\frac{s}{c_{2}%
}\right)  \int_{\Omega}\left(  -\eta_{i}\right)  \left\vert f_{i}\right\vert
^{2}\text{.}
\]

Thanks again to Lemma \ref{lemtheta} $\left(  ii\right)  -\left(  iii\right)
$, the properties of $\varphi_{3}$ require to treat separately the cases when
$x\in\vartheta$ and $x\in\Omega\left\backslash \vartheta\right.  $, where
$\vartheta\subset\Omega$ is a neighborhood of $\partial\Omega$ given by
$\vartheta=\left\{  \rho\leq\left\vert x\right\vert <R\right\}  $, with
$\rho>0$ such that $x_{0}\notin\overline{\vartheta}$. We take here $\rho
=\frac{|x_{0}|+R}{2}$. We first observe that when $x\in\vartheta$ ,
\begin{equation}
-\left\vert \nabla\Phi_{3}\left(  x,\cdot\right)  \right\vert ^{2}%
=-\frac{s^{2}}{\Gamma^{2}}\left\vert \nabla\varphi_{3}\left(  x\right)
\right\vert ^{2}\leq-\frac{s^{2}}{c_{2}\Gamma^{2}}\left\vert \varphi
_{3}\left(  x\right)  \right\vert \leq\frac{2s}{c_{2}}\eta_{3}\left(
x\right)  \text{ .} \label{nez}%
\end{equation}

We see that
\[%
\begin{array}
[c]{ll}%
\displaystyle\sum_{i=3,4}\left(  -\frac{1}{8}d_{i}\displaystyle\int_{\Omega
}\left\vert \nabla\Phi_{i}\right\vert ^{2}\left\vert f_{i}\right\vert
^{2}\right)  & \leq\displaystyle\sum_{i=3,4}\left(  -\frac{1}{8}%
d_{i}\displaystyle\int_{\vartheta}\left\vert \nabla\Phi_{i}\right\vert
^{2}\left\vert f_{i}\right\vert ^{2}\right)  \text{ since }\vartheta
\subset\Omega\\
& \leq\displaystyle\sum_{i=3,4}\frac{1}{4}d_{i}\frac{s}{c_{2}}%
\displaystyle\int_{\vartheta}\eta_{i}\left\vert f_{i}\right\vert ^{2}\text{
thanks to estimate (\ref{nez})}\\
& =\displaystyle\sum_{i=3,4}\frac{1}{4}d_{i}\frac{s}{c_{2}}\displaystyle\int
_{\Omega}\eta_{i}\left\vert f_{i}\right\vert ^{2}-\displaystyle\sum
_{i=3,4}\frac{1}{4}d_{i}\frac{s}{c_{2}}\displaystyle\int_{\Omega
\left\backslash \vartheta\right.  }\eta_{i}\left\vert f_{i}\right\vert ^{2}\\
& \leq\displaystyle\sum_{i=3,4}\left(  -\frac{1}{4}d_{i}\frac{s}{c_{2}%
}\right)  \displaystyle\int_{\Omega}\left(  -\eta_{i}\right)  \left\vert
f_{i}\right\vert ^{2}+C_{7}\displaystyle\sum_{i=1,2}\displaystyle\int_{\Omega
}\left\vert f_{i}\right\vert ^{2}%
\end{array}
\]
where in the last line, we defined
\[
C_{6}:=\underset{\overline{\Omega}}{\text{max}}\left\vert \psi\right\vert
+\underset{i}{\text{max}}\,d_{i}\, \underset{\overline{\Omega}}{\text{max}%
}\left\vert \nabla\psi\right\vert ^{2}\text{ ,}\qquad C_{7}:=\underset
{i}{\text{max}}\,d_{i}\,\frac{C_{6}}{c_{2}\,c_{3}^{2}}\text{ ,}%
\]
and noticed that%
\[%
\begin{array}
[c]{ll}%
\displaystyle-\frac{1}{4}d_{3}\frac{s}{c_{2}}\int_{\Omega\left\backslash
\vartheta\right.  }\eta_{3}\left\vert f_{3}\right\vert ^{2} & \leq
\displaystyle\frac{1}{4}d_{3}\frac{s}{c_{2}}\int_{\Omega\left\backslash
\vartheta\right.  }\left(  \frac{C_{6}s}{\Gamma^{2}}\right)  \left\vert
u_{3}\right\vert ^{2}e^{s\frac{1}{\Gamma}\varphi_{3}}\text{ since }\left\vert
\eta_{3}\right\vert \leq\frac{C_{6}s}{\Gamma^{2}}\\
& =\displaystyle\frac{1}{4}d_{3}\frac{s}{c_{2}}\int_{\Omega\left\backslash
\vartheta\right.  }\left(  \frac{C_{6}s}{\Gamma^{2}}\right)  \left\vert
u_{3}\right\vert ^{2}e^{s\frac{1}{\Gamma}\varphi_{1}}e^{s\frac{1}{\Gamma
}\left(  \varphi_{3}-\varphi_{1}\right)  }\\
& \leq\displaystyle\frac{1}{4}d_{3}\frac{s}{c_{2}}\int_{\Omega\left\backslash
\vartheta\right.  }\left(  \frac{C_{6}s}{\Gamma^{2}}\right)  \left\vert
u_{3}\right\vert ^{2}e^{s\frac{1}{\Gamma}\varphi_{1}}e^{-s\frac{1}{\Gamma
}c_{3}}\text{ by Lemma \ref{lemtheta}}(iii)\\
& \leq C_{7}\displaystyle\int_{\Omega\left\backslash \vartheta\right.
}\left\vert u_{3}\right\vert ^{2}e^{s\frac{1}{\Gamma}\varphi_{1}}%
=C_{7}\displaystyle\int_{\Omega\left\backslash \vartheta\right.  }\left\vert
u_{1}\right\vert ^{2}e^{s\frac{1}{\Gamma}\varphi_{1}}\text{ because }%
u_{3}=u_{1}\\
& \leq C_{7}\displaystyle\int_{\Omega}\left\vert f_{1}\right\vert ^{2}\text{
.}%
\end{array}
\]
We proceed similarly for $i=4$ and get $\displaystyle-\frac{1}{4}d_{4}\frac
{s}{c_{2}}\int_{\Omega\left\backslash \vartheta\right.  }\eta_{4}\left\vert
f_{4}\right\vert ^{2}\leq C_{7}\displaystyle\int_{\Omega}\left\vert
f_{2}\right\vert ^{2}$. Finally, we see that
\[
\sum_{i=1,..,4}\left(  -\frac{1}{8}d_{i}\int_{\Omega}\left\vert \nabla\Phi
_{i}\right\vert ^{2}\left\vert f_{i}\right\vert ^{2}\right)  \leq
\sum_{i=1,..,4}\left(  -\frac{1}{4}d_{i}\frac{s}{c_{2}}\right)  \int_{\Omega
}\left(  -\eta_{i}\right)  \left\vert f_{i}\right\vert ^{2}+C_{7}\sum
_{i=1,2}\int_{\Omega}\left\vert f_{i}\right\vert ^{2}\text{ .}%
\]
The fact that $\frac{1}{\Gamma}\leq\frac{1}{h}$ completes the proof of Lemma
\ref{prop4}. \medskip

Consequently, by (\ref{4.1}), (\ref{4.2}), Lemma \ref{lemnn} and Lemma
\ref{prop4}, for any $s\in\left(  0,s_{1}\right]  \cap\left(  0,s_{0}\right]
$ and any $h\in\left(  0,1\right]  $, we see that
\[%
\begin{array}
[c]{ll}%
\left\langle \mathcal{S}^{\prime}f,f\right\rangle +2\left\langle
\mathcal{S}f,\mathcal{A}f\right\rangle  & \leq(C_{3}+C_{5})\displaystyle\frac
{s}{\Gamma}\displaystyle\sum_{i=1,..,4}d_{i}\displaystyle\int_{\Omega}|\nabla
f_{i}|^{2}+(C_{3}+C_{5}+C_{7})\displaystyle\frac{1}{h^{2}}\left\Vert
f\right\Vert ^{2}\\
& \quad+\displaystyle\frac{1}{\Gamma}\displaystyle\sum_{i=1,..,4}\left(
2-\displaystyle\frac{d_{i}s}{4c_{2}}\right)  \displaystyle\int_{\Omega}%
(-\eta_{i})\left\vert f_{i}\right\vert ^{2}\\
& \leq(C_{3}+C_{5})\displaystyle\frac{s}{\Gamma}\displaystyle\sum
_{i=1,..,4}d_{i}\displaystyle\int_{\Omega}|\nabla f_{i}|^{2}+(C_{3}%
+C_{5}+C_{7})\displaystyle\frac{1}{h^{2}}\left\Vert f\right\Vert ^{2}\\
& \quad+\left(  2-\underset{i}{\text{min}}\,d_{i}\displaystyle\frac{s}{4c_{2}%
}\right)  \displaystyle\frac{1}{\Gamma}\displaystyle\sum_{i=1,..,4}%
\displaystyle\int_{\Omega}(-\eta_{i})\left\vert f_{i}\right\vert ^{2}\text{ .}%
\end{array}
\]
For $s\in(0,s_{2}]$, where $s_{2}:=$min$(s_{0},s_{1},(C_{3}+C_{5})^{-1}%
,c_{2}\,(\underset{i}{\text{min}}\,d_{i})^{-1})$, we see that
\[
\left\langle \mathcal{S}^{\prime}f,f\right\rangle +2\left\langle
\mathcal{S}f,\mathcal{A}f\right\rangle \leq\frac{1+C_{0}}{\Gamma}\left\langle
\mathcal{S}f,f\right\rangle +\frac{C_{1}}{h^{2}}\left\Vert f\right\Vert
^{2}\text{ ,}%
\]
with $C_{0}:=1-\underset{i}{\text{min}}\,d_{i}\,\displaystyle\frac{s_{2}}%
{4c_{2}}\in\left(  0,1\right)  $ and $C_{1}\geq C_{3}+C_{5}+C_{7}$. We
complete the proof of Proposition~\ref{prop2} by taking
\[
C_{1}:=\text{max}(1,C_{3}+C_{5}+C_{7},4K_{0}\,(1+K_{0}\,C_{Sob}),4K_{0}%
^{2}\,C_{Sob}(\underset{i}{\text{min}}\,d_{i})^{-1})\text{ .}%
\]

\bigskip

\bigskip

$\bigskip$

\bigskip

\bigskip

\bigskip

\bigskip

\bigskip

\bigskip

\begin{thebibliography}{9999}                                                                                             %


\bibitem[BP]{BP}C. Bardos and K.D. Phung, Observation estimate for kinetic
transport equations by diffusion approximation. C. R. Math. Acad. Sci. Paris,
355, no.6, (2017), 640--664.

\bibitem[BT]{BT}C. Bardos and L. Tartar, Sur l'unicit\'{e} retrograde des
\'{e}quations paraboliques et quelques questions voisines. Arch. Rational
Mech. Anal., 50 (1973), 10--25.

\bibitem[BDS]{BDS}M. Bisi, L. Desvillettes and G. Spiga, Exponential
Convergence to Equilibrium via Lyapounov Functionals for Reaction-Diffusion
Equations Arising from non Reversible Chemical Kinetics. Mathematical
Modelling and Numerical Analysis, vol. 43, n.1, (2009), 151--172.

\bibitem[BuP]{BuP}R. Buffe and K.D. Phung, Observation estimate for the heat
equations with Neumann boundary condition via logarithmic convexity. arXiv:2105.12977

\bibitem[D]{D}L. Desvillettes, About Entropy Methods for Reaction-Diffusion
Equations. Rivista di Matematica dell'Universit\`{a} di Parma, vol. 7, n.7,
(2007), 81--123.

\bibitem[DF]{DF}L. Desvillettes and K. Fellner, Exponential decay toward
equilibrium via entropy methods for reaction-diffusion equations. Journal of
Mathematical Analysis and Applications, vol. 319, n.1, (2006), 157--176.

\bibitem[DF2]{DF2}L. Desvillettes and K. Fellner, Entropy Methods for
Reaction-Diffusion Systems. Discrete and Continuous Dynamical Systems,
supplement (2007), 304--312.

\bibitem[DF3]{DF3}L. Desvillettes and K. Fellner, Entropy Methods for
Reaction-Diffusion Equations: Slowly Growing A-priori Bounds. Revista
Matematica Iberoamericana, vol. 24, n.2, (2008), 407--431.

\bibitem[DFT]{DFT}L. Desvillettes, K. Fellner and B. Q. Tang, Trend to
equilibrium for reaction-diffusion systems arising from complex balanced
chemical reaction networks. SIAM Journal of Mathematical Analysis, vol. 49,
(2017), 2666--2709.

\bibitem[HZ]{HZ}V. Hern\'{a}ndez-Santamar\'{\i}a and E. Zuazua,
Controllability of shadow reaction-diffusion systems. J. Differ. Equ. 268
(2020), 3781--3818.

\bibitem[LSU]{LSU}O.A. Ladyzenskaja, V.A. Solonnikov, and N.N. Uralceva,
\textit{Linear and quasilinear equations of parabolic type}. Volume 23.
American Mathematical Soc., 1968.

\bibitem[Le B]{Le B}K. Le Balc'h, Controllability of a 4x4 quadratic
reaction-diffusion system. J. Differ. Equ. 266(6) (2019), 3100--3188.

\bibitem[P]{P}K.D. Phung, Carleman commutator approach in logarithmic
convexity for parabolic equations. Mathematical Control and Related Fields 8
(3-4) (2018), 899--933.

\bibitem[PW]{PW}K.D. Phung and G. Wang, An observability estimate for
parabolic equations from a measurable set in time and its applications.
Journal of the European Mathematical Society 15 (2) (2013) 681--703.
\end{thebibliography}
\end{document}